\documentclass[10pt]{article}

\usepackage{ucs}
\usepackage[utf8x]{inputenc}
\usepackage{graphicx}
\usepackage{amsfonts}
\usepackage{dsfont}
\usepackage{amssymb}
\usepackage{amsmath}
\usepackage{amsthm}
\usepackage{enumerate}
\usepackage{stmaryrd}
\usepackage{fullpage}
\usepackage{ifthen}
\usepackage{subfigure}
\usepackage{epic}
\usepackage{authblk}
\usepackage{textcomp}
\usepackage[all]{xy}
\usepackage{mathrsfs}

\usepackage[hypertexnames=false,colorlinks=true,linkcolor=blue,citecolor=blue]{hyperref}
\usepackage[numbers,comma,square,sort&compress]{natbib}
\usepackage[letterpaper,margin=1in,centering]{geometry}

\usepackage{color}
\graphicspath{{eps/}{pdf/}}

\newcommand{\bqq}{\begin{equation}}
\newcommand{\eqq}{\end{equation}}

\newcommand{\C}{\mathbb{C}}

\newcommand{\N}{\mathbb{N}}
\newcommand{\R}{\mathbb{R}} 
\newcommand{\Z}{\mathbb{Z}}

\newcommand{\B}{\mathcal{B}}
\newcommand{\cC}{\mathscr{C}}
\newcommand{\cD}{\mathcal{D}}
\newcommand{\E}{\mathcal{E}}
\newcommand{\F}{\mathcal{F}}
\newcommand{\G}{\mathcal{G}}

\newcommand{\I}{\mathcal{I}}
\newcommand{\J}{\mathcal{J}}
\newcommand{\K}{\mathcal{K}}
\newcommand{\cL}{\mathcal{L}}
\newcommand{\M}{\mathcal{M}}
\newcommand{\cN}{\mathcal{N}}
\newcommand{\cO}{\mathcal{O}}
\newcommand{\co}{\mathrm{o}}

\newcommand{\Q}{\mathcal{Q}}
\newcommand{\cR}{\mathcal{R}}
\newcommand{\cS}{\mathcal{S}}
\newcommand{\T}{\mathcal{T}}
\newcommand{\U}{\mathcal{U}}
\newcommand{\V}{\mathcal{V}}

\newcommand{\X}{\mathcal{X}}
\newcommand{\Y}{\mathcal{Y}}
\newcommand{\cZ}{\mathcal{Z}}

\newcommand{\vp}{\varphi}

\newcommand{\rK}{\mathscr{K}}
\newcommand{\rG}{\mathscr{G}}

\newcommand{\mS}{\mathbf{S}}

\newcommand{\bA}{\overline{A}}
\newcommand{\bB}{\overline{B}}
\newcommand{\bC}{\overline{C}}
\newcommand{\bD}{\overline{D}}
\newcommand{\bE}{\overline{E}}
\newcommand{\bF}{\overline{F}}

\newcommand{\bu}{\mathbf{u}}
\newcommand{\bv}{\mathbf{v}}
\newcommand{\be}{\mathbf{e}}
\newcommand{\md}{\mathrm{d}}
\newcommand{\rme}{\mathrm{e}}
\newcommand{\rmi}{\mathrm{i}}
\newcommand{\mf}{\mathbf{f}}

\numberwithin{equation}{section}

\newtheorem{lem}{Lemma}[section]
\newtheorem{thm}{Theorem}

\newtheorem{cor}{Corollary}
\newtheorem{rmk}[lem]{Remark}

\newtheorem{hyp}[lem]{Hypothesis}

\newenvironment{Hypothesis}[1]%
  {\begin{trivlist}\item[]{\bf Hypothesis #1 }\em}{\end{trivlist}}

\makeatletter\@addtoreset{figure}{section}\makeatother

\makeatletter \@addtoreset{equation}{section} \makeatother

\newenvironment{Proof}[1][.]%
 {\begin{trivlist}\item[]\textbf{Proof#1 }}%
 {\hspace*{\fill}$\rule{0.3\baselineskip}{0.35\baselineskip}$\end{trivlist}}

\title{Center Manifolds without a Phase Space }
\author[1]{Gr\'egory Faye\footnote{Corresponding author, gregory.faye@math.univ-toulouse.fr}}
\affil[1]{\small CNRS, UMR 5219, Institut de Math\'ematiques de Toulouse, 31062 Toulouse Cedex, France}
\author[2]{Arnd Scheel\footnote{AS was partially supported by the National Science Foundation through grants NSF-DMS-1311740 and NSF-DMS-1612441, and through a DAAD Fellowship.}}
\affil[2]{\small University of Minnesota, School of Mathematics, 206 Church Street S.E., Minneapolis, MN 55455, USA}

\begin{document}
\maketitle

\begin{abstract}
We establish center manifold theorems that allow one to study the bifurcation of small solutions from a trivial state in systems of functional equations posed on the real line. The class of equations includes most importantly  nonlinear equations with nonlocal coupling through convolution operators as they arise in the description of spatially extended dynamics in neuroscience. These systems possess a natural spatial translation symmetry but local existence or uniqueness theorems for a spatial evolution associated with this spatial shift or even a well motivated choice of phase space for the induced dynamics do not seem to be available, due to the infinite range forward- and backward-coupling through nonlocal convolution operators. We perform a reduction relying entirely on  functional analytic methods.  Despite the nonlocal nature of the problem, we do recover a local differential equation describing the dynamics on the set of small bounded solutions, exploiting that the translation invariance of the original problem induces a flow action on the center manifold. We apply our reduction procedure to problems in mathematical neuroscience, illustrating in particular the new type of algebra necessary for the computation of Taylor jets of reduced vector fields.

\end{abstract}


%\tableofcontents

\section{Introduction}
Center-manifold reductions have become a central tool to the analysis of dynamical systems. The very first results on center manifolds go back to the pioneering works of Pliss \cite{pliss:64} and Kelley \cite{kelley:67} in the finite-dimensional setting. In the simplest context, one studies differential equations in the vicinity of a non-hyperbolic equilibrium, 
\[
\frac{\md u}{\md t}=f(u)\in\R^n, \qquad f(0)=0, \quad \mathrm{spec}(f'(0))\cap \rmi\R\neq \emptyset.
\]
The basic reduction establishes that the set of small bounded solutions $u(t)$, $t\in\R$, $\sup |u(t)|<\delta\ll 1$, is \emph{pointwise} contained in a manifold, that is,  $u(t)\in W^\mathrm{c}$ for all $t$. This manifold is a subset of phase space, $W^\mathrm{c}\subset\R^n$, contains the origin, $0\in W^\mathrm{c}$, and is tangent to $E^\mathrm{c}$, the generalized eigenspace associated with purely imaginary eigenvalues of $f'(0)$. As a consequence, the flow on $W^\mathrm{c}$ can be projected onto $E^\mathrm{c}$, to yield a \emph{reduced vector field}. The reduction to this lower-dimensional ODE then allows one to describe solutions qualitatively, even explicitly in some cases. Of course, the method applies to higher-order differential equation, which one simply writes as first-order equation in a canonical fashion.  Extensions to infinite-dimensional dynamical systems were pursued soon after; see for instance \cite{henry:81}. 

Starting with the work of Kirchg\"assner \cite{ki82}, such reductions have been extended to systems with $u\in \X$, a Banach space, where the initial value problem is not well-posed: For most initial conditions $u_0$, there does not exist a local solution $u(t)$, $0\leq t<\delta$, say. Local solutions do exist however for all initial conditions on a finite-dimensional center-manifold, and much of the theory is quite analogous to the finite-dimensional case; see \cite{vander-iooss:92}. In these theories, one can typically split the phase space in infinite-dimensional linear spaces where solutions to the linearized equation either decay or grow, and a finite-dimensional center subspace. Such splittings are known as Wiener-Hopf factorizations and can be  difficult to achieve in the case of forward-backward delay equations, where nevertheless center-manifold reductions are available \cite{hv}. 

Our point of view here is slightly more abstract, shedding the concept of a phase space in favor of a focus on small bounded trajectories. We perform a purely functional analytic reduction, based on Fredholm theory \cite{faye-scheel:13} in the space of bounded trajectories (rather than the phase space). We parameterize the set of bounded solutions by the set of (weakly) bounded solutions to the linear equation, which is a finite-dimensional vector space, amenable to a variety of parameterizations. Only after this reduction, we derive a differential equation on this finite-dimensional vector space, whose solutions, when lifted to the set of bounded solutions to the nonlinear problem describe all small bounded solutions. 

To be more precise, we focus on nonlocal equations of the form 
\bqq
 u+\K \ast u +\F(u)=0,
\label{eqNL}
\eqq
for $u:\R\rightarrow\R^n$, $n\geq1$. Here, $\K\ast u$ stands for matrix convolution on $\R$,
\begin{equation*}
(\K\ast u(x))_i = \sum_{j=1}^n  \int_\R \K_{i,j}(x-y)u_j(y)\md y,\quad 1\leq i\leq n,
\end{equation*}
and $\F(u)$ encodes nonlinear terms, possibly also involving nonlocal interactions.  

A prototypical example arises when studying stationary or traveling-wave solutions to neural field equations, which are used in mathematical neuroscience to model cortical activity. A typical model is
\bqq
\frac{\md u}{\md t} = - u+\K \ast S(u),
\label{nfe}
\eqq
where $u(t,x)\in\R$ represents a locally  averaged membrane potential,  the nonlinearity $S$ denotes a firing rate function and the kernel $\K$ encodes the connectivity, \textit{i.e.} how neurons located at position $x$ interact with neurons located at position $y$ across the cortex. Stationary solutions of \eqref{nfe} are thought to be associated to short term memory and to encode our ability to remember given tasks over a period of milliseconds, providing motivation for extensive studies of such solutions over the past two decades; see for instance  \cite{faye-etal:13} for a more thorough presentation of the problem and  related references.

Beyond techniques based on comparison principles, which apply to some extent when, say, $\K>0$, a widely used method to study the stationary problem
\bqq
0 = - u+\K \ast S(u),
\label{nfes}
\eqq
focuses on kernels $\K$ with rational Fourier transform, 
\begin{equation*}
\widehat{\K}(\rmi \ell) = \int_\R \K(x) \rme^{-\rmi \ell x} \md x = \frac{\mathrm{Q}(\ell^2)}{\mathrm{P}(\ell^2)},
\end{equation*}
for some polynomials $\mathrm{P}$ and $\mathrm{Q}$ with $\text{deg}\,\mathrm{Q} < \text{deg}\,\mathrm{P}$; see \cite{laing-troy:03}. 
There appears to be little motivation for such special kernels other than the obvious technical advantage that the nonlocal equation can be written as a local differential equation,
\bqq
0 = - \mathrm{Q}(-\partial_{xx})u+\mathrm{P}(-\partial_{xx}) S(u),
\label{nfesq}
\eqq
where now dynamical systems techniques, in particular center manifold reduction, are applicable. On the other hand, 
$\K>0$ restricts to the very specific case of excitatory connections, and scalar (or, say, cooperative) dynamics.

While our results are motivated to some extent by the desire to eliminate the unnecessary restriction to rational Fourier transforms, we believe that there is more generally valuable insight in the results presented here. For instance, many problems with nonlocal, pseudo-differential operators can be cast in the form \eqref{nfes}, after possibly preconditioning the equation with the resolvent of a leading-order part.

\paragraph{Summary of main results.}

We now state our main result in a somewhat informal way. We study 
\bqq
\T u +\F(u)=0,\qquad \T u = u+\K\ast u,
\label{eqNLtime}
\eqq
\begin{itemize}
\item \emph{Exponential Localization}: the interaction kernel $\K$ and its derivative $\K'$ are exponentially localized (see Section \ref{subsec:hyp}, Hypothesis (H1));
\item \emph{Smoothness and Invariance}: the nonlinear operator $\F$ is assumed to be sufficiently smooth and translation invariant $\F(u(\cdot + \xi))(\cdot)=\F(u(\cdot ))(\cdot+\xi)$, with $\F(0)=0$, $D_u\F(0)= 0$ (see Section \ref{subsec:hyp} Hypothesis (H2)).
\end{itemize}
Using Fourier transform, one can readily find the finite-dimensional space $\ker \T$ of solutions to $\T u=0$ with at most algebraic growth and construct a bounded projection $\Q$ onto this set, in a space of functions allowing for slow exponential growth. 

\begin{thm}\label{thmCM}
Assume that the interaction kernel $\K$ and the nonlinear operator $\F$ satisfy Hypotheses (H1)-(H2). Then, there exists $\delta>0$  and a map $\Psi\in \cC^k(\ker\T,\ker \Q)$ with $\Psi(0)=D_u\Psi(0)=0$,  such that the manifold
\begin{equation*}
\M_0:=\left\{ u_0+\Psi(u_0) ~\middle|~ u_0 \in \ker \T \right\} 
\end{equation*}
contains the set of all bounded solutions of \eqref{eqNLtime} with  $\sup_{x\in\R}| u(x) |\leq \delta$. 
\end{thm}
We refer to $\M_0$ as a (global) center manifold for \eqref{eqNLtime}. Note however that points on $\M_0$ consist of \emph{trajectories}, that is, of solutions $u(x)$, $x\in\R$, rather than of initial values to solutions, in the more common view of center manifolds. Also note that, according to the theorem, $\M_0$ only contains the set of bounded solutions, not all elements of $\M_0$ are necessarily bounded solutions. As is well known from the classical center manifold theorem, the set of bounded solutions may well be trivial, consisting of the point $u\equiv 0$, only, rather than being diffeomorphic to a finite-dimensional ball. It is therefore necessary to study the elements of $\M_0$ in more detail. 

We will see in the proof that, as is common in the construction of center manifolds, we modify the nonlinearity $\F$ to $\F_\epsilon$ outside of a small $\epsilon$-neighborhood, $\sup_x|u(x)|\leq \epsilon$, in the construction of $\M_0$. Therefore, \emph{all} elements of $\M_0$ are in fact solutions, with possibly mild exponential growth, and to the modified equation $\T u + \F_\epsilon(u)=0$. The set of solutions to this equation is translation invariant and can be described by a differential equation as stated in the following result.

\begin{cor}\label{corRE}
Under the assumptions (H1)-(H2) of Theorem \ref{thmCM}, any element  $u=u_0+\Psi(u_0)$ of $\M_0$, corresponds to a unique solution of a differential equation
\bqq
\label{eqReduced}
\frac{\md u_0}{\md x}=f(u_0),
\eqq
on the linear vector space $u_0\in \ker\T$. 
The Taylor jet of $f$ can be computed from properties of $\T$ and $\F$, solving linear equations, only. 
\end{cor}

Note that the differentiation in \eqref{eqReduced} does not refer to differentiation of $u_0$, which of course is a function of $x$ when viewed as an element of the kernel. We rather view $\ker\T$ as an abstract vector space on which we study the differential equation \eqref{eqReduced}. Note also that we do not claim that every solution to \eqref{eqReduced} is a solution to \eqref{eqNLtime} --- this is true only for small solutions.

We will explain below how to actually compute the Taylor jet of $f$. Having access to \eqref{eqReduced} as a means of describing elements of $\M_0$, the abstract reduction Theorem \ref{thmCM} becomes very valuable: one simply studies the differential equation \eqref{eqReduced}, or, to start with, the equation obtained from the leading order Taylor approximation,  using traditional dynamical systems methods.  Small bounded solutions obtained in this fashion will then correspond to solutions of the original nonlocal problem \eqref{eqNLtime}.

\paragraph{Outline.} We state a precise version of  our main theorem on the existence of a center manifold for systems of nonlocal equations, mention extensions, and provide basic tools necessary for the application  in Section \ref{sec:main}. Proofs are given in Section \ref{s:proofs}, applications to neural field equations in Section \ref{sec:appli}.

\section{Existence of center manifolds ---  main result and extensions}\label{sec:main}

We introduce the functional analytic framework, and state the main hypotheses on linear and nonlinear parts of the equation in Section \ref{subsec:hyp}. We then state the main theorem of this paper, Section \ref{s:main}, and extensions, Section \ref{s:ext}. Sections \ref{s:non}--\ref{subsec:example} provide basic tools that allow one to apply the main result, showing how to verify assumptions on the nonlinearity, Section \ref{s:non}, how to construct projections $\Q$, Section \ref{s:pro}, and how to compute Taylor jets of the reduced vector field, Section \ref{subsec:example}. 

\subsection{Functional-analytic setup and main assumptions}\label{subsec:hyp}

We introduce function spaces and state our main hypotheses, (H1) and (H2). 
\paragraph{Function spaces.}
For $\eta \in \R$, $1\leq p\leq\infty$, we define the weighted space $L^p_\eta(\R,\R^n)$, or simply $L^p_\eta$, when $n=1$, through
\begin{equation*}
L^p_\eta(\R,\R^n):= \left\{u\in L^p_\mathrm{loc}(\R,\R^n) ~:~ \omega_\eta u \in L^p\left(\R,\R^n\right)\right\},
\end{equation*}
where $\omega_\eta$ is a $\cC^\infty$ function defined as
\begin{equation*}
\omega_\eta (x) =\left\{\begin{array}{lcl}
\rme^{\eta x} & \text{for} & x \geq 1,\\
\rme^{-\eta x} & \text{for} & x \leq -1
\end{array}
 \right., \quad \omega_\eta>0 \text{  on } [-1,1].
\end{equation*}
We also use the standard Sobolev spaces $W^{k,p}(\R,\R^n)$, or simply $W^{k,p}$ when $n=1$, for $k \geq 0$ and $1\leq p \leq \infty$:
\begin{equation*}
W^{k,p}(\R,\R^n):=\left\{u \in L^p\left(\R,\R^n\right)~:~ \partial^\alpha_x u \in L^p\left(\R,\R^n\right), \quad 1 \leq \alpha \leq k \right\},
\end{equation*}
with norm
\begin{equation*}
\| u \|_{W^{k,p}(\R,\R^n)}=\left\{
\begin{array}{cc}
\left( \sum_{\alpha\leq k}  \| \partial^\alpha_x u \|_{L^p(\R,\R^n)}^p \right)^{\frac{1}{p}}, & 1 \leq p <\infty \\
\underset{\alpha \leq k}{\max}  \| \partial^\alpha_x u \|_{L^\infty(\R,\R^n)}, & p = \infty.
\end{array}
\right.
\end{equation*}
We denote by $H^k(\R,\R^n)$ the Sobolev space $W^{k,2}(\R,\R^n)$ and use the weighted spaces  $W^{k,p}_\eta(\R,\R^n)$ and $H^k_\eta(\R,\R^n)$ the weighted Sobolev spaces defined through
\begin{equation*}
W^{k,p}_\eta(\R,\R^n):=\left\{u \in L^p_\mathrm{loc}\left(\R,\R^n\right)~:~ \omega_\eta \partial^\alpha_x u \in L^p\left(\R,\R^n\right), \quad 1 \leq \alpha \leq k \right\},
\end{equation*}
with norm
\begin{equation*}
\| u \|_{W^{k,p}_\eta(\R,\R^n)}=\left\{
\begin{array}{cc}
\left( \sum_{\alpha\leq k}  \| \omega_\eta \partial^\alpha_x u \|_{L^p(\R,\R^n)}^p \right)^{\frac{1}{p}}, & 1 \leq p <\infty \\
\underset{\alpha \leq k}{\max}  \| \omega_\eta \partial^\alpha_x u \|_{L^\infty(\R,\R^n)}, & p = \infty,
\end{array}
\right.
\end{equation*}
and $H^k_\eta(\R,\R^n):=W^{k,2}_\eta(\R,\R^n)$.

\paragraph{Assumptions on the linear part.}

We require that the convolution kernel is exponentially localized and smooth in the following sense. 
\begin{Hypothesis}{(H1)}
We assume that there exists $\eta_0>0$ such that $\K_{i,j} \in W^{1,1}_{\eta_0}(\R)$ for all $1\leq i,j \leq n$.
\end{Hypothesis}
We define the complex Fourier transform $\widehat{\K}(\nu)$ of $\K$ as
\bqq
\label{Fourier}
\widehat{\K}(\nu)=\int_\R \K(x) \rme^{-\nu x} \md x,
\eqq
for all $\nu\in\C$ where the above integral is well-defined. Note that because each component of the matrix kernel $\K $ belongs to $L^1_{\eta_0}$, the Fourier transform $\widehat{\K}(\nu)$ is analytic in the strip $\mathscr{S}_{\eta_0}:=\left\{ \nu\in\C ~|~ |\Re(\nu)|<\eta_0 \right\}$. Thereby, the characteristic equation
\bqq
\label{eqDisp}
d(\nu):=\det\left(I_n+\widehat{\K}(\nu)\right)=0
\eqq
is an analytic function in the strip and has isolated roots on the imaginary axis, when counted with multiplicity. Moreover, since  $\K' \in L^1_{\eta_0}$ (component-wise), we have $|\widehat{\K}(\rmi \ell+\eta)|\underset{\ell \rightarrow \pm\infty}{\longrightarrow}0 $, for $|\eta|<\eta_0$, such that the number of roots of $d$ on the imaginary axis, counted with multiplicity, is finite. Throughout, we will assume that the number of roots is not zero, in which case our results would be trivial. 
 
We consider $\T$ as a bounded operator on $H^1_{-\eta}(\R,\R^n)$, $0<\eta<\eta_0$, slightly abusing notation and not making the dependence of $\T$ on $\eta$ explicit. With the natural bounded inclusion  $\iota^{\eta,\eta'}$, $\eta<\eta'$, one finds $\T \iota^{\eta,\eta'}=\iota^{\eta,\eta'} \T$.  Now, by finiteness of the number of roots of $d$, we can choose $\eta_1>0$, small, such that $d(\nu)$ does not vanish in $0<|\Im\nu|\leq \eta_1$. We will then find that the kernel $\E_0$ of $\T$ is independent of $\eta$ for $0<\eta<\eta_1$ in the sense that $\iota^{\eta,\eta'}$ provides isomorphisms between kernels for $\eta$ and $\eta'$, $0<\eta<\eta'<\eta_1$. The dimension of $\E_0$ is given by the sum of multiplicities of roots $\nu\in\rmi \R$ of $d(\nu)$, with a basis of the form $p(x)\rme^{\nu x}$,  $p$ a vector-valued polynomial of degree at most $m-1$  when $\nu$ is a root of $d$ of order $m$; see Lemma \ref{l:21}, below.   We also need a bounded projection 
\begin{equation}\label{e:pro}
\Q:H^1_{-\eta}(\R,\R^n)\rightarrow H^1_{-\eta}(\R,\R^n), \qquad \Q^2=\Q,\qquad  \mathrm{rg}\,(\Q)=\E_0=\ker\T,
\end{equation}
with a continuous extension to $L^2_{-\eta}(\R,\R^n)$. Again, we require  $\Q \iota^{\eta,\eta'} = \iota^{\eta,\eta'}  \Q$, a possible choice being the $L^2_{\eta_1}(\R,\R^n)$-orthonormal projection.  We discuss this in more detail in Section \ref{s:pro}, including computationally advantageous choices of projections.

\paragraph{Assumptions on the nonlinear part.}

A common approach to the construction of center manifolds is to modify the nonlinearity outside of a small neighborhood of the origin. We therefore first define a pointwise, smooth cut-off function $\bar{\chi}:{\R^n}\to\R$, with
\[
\bar{\chi}(u)= \left\{\begin{array}{lcl} 1 & \text{for} & \| u \|\leq 1 \\ 0 & \text{for} & \| u \| \geq 2 \end{array} \right.,\quad \bar{\chi}(u)\in [0,1],
\]
and then a cut-off operator  $\chi_\epsilon$, mapping measurable functions $u:\R\to \R^n$ into $L^\infty(\R,\R^n)$,
\[
{\chi_\epsilon(u)(x)=\bar{\chi}(u(x)/\epsilon)\cdot u(x).}
\]
Lastly, formally define the family of translation operators $\tau_\xi$, $\xi\in\R$,
\[
(\tau_\xi\cdot u)(x):=u(x-\xi),
\]
the canonical representation of the group $\R$ on functions over $\R$. Slightly abusing notation, we will use the same symbol $\tau_\xi$ for the action on various function spaces. Note that $\tau_\xi$ will be bounded for $
\xi$ fixed on all spaces introduced above. We define the modified nonlinearities 
\begin{equation}\label{modF}
\F^\varepsilon:=\F\circ \chi_\epsilon.
\end{equation}

\begin{Hypothesis}{(H2)}
We assume that there exists $k\geq 2$ and $\eta_0>0$ such that for all $\epsilon>0$, sufficiently small, the following properties hold.
\begin{enumerate}
\item $\F\in \cC^k(\V,W^{1,\infty}(\R,\R^n)$, for some small neighborhood $\V\subset W^{1,\infty}(\R,\R^n)$, and $\F(0)=0$, $D_u\F(0)= 0$;
\item $\F$ commutes with translations, $\F\circ\tau_\xi=\tau_\xi\circ\F$ for all $\xi\in\R$;
\item $\F^\epsilon:  H^1_{-\zeta}(\R,\R^n) \longrightarrow H^{1}_{-\eta }(\R,\R^n)$ is $\cC^k$ for all nonnegative pairs $(\zeta,\eta)$ such that $0<k\zeta <\eta<\eta_0$,  $D^j\F^\epsilon(u): (H^1_{-\zeta}(\R,\R^n))^j \longrightarrow H^{1}_{-\eta }(\R,\R^n)$ is bounded for $0<j\zeta\leq \eta<\eta_0$, $0\leq j\leq k$ and Lipschitz in $u$ for $1\leq j \leq  k-1$.
\end{enumerate}
\end{Hypothesis}
Note that $\F^\epsilon$  commutes with $\tau_\xi$ since $\F$ and $\chi_\epsilon$ do. The first condition is the common condition, guaranteeing that $\T$ is actually the linearization at an equilibrium $u\equiv 0$, that is, at a solution invariant under translations $\tau_\xi$. The second condition puts us in the scenario of an autonomous dynamical system. 
The last condition on the modified nonlinearity is a technical condition, known from the proofs of smoothness of center manifolds in ODEs \cite{vander-vangils:87}, that will imply smoothness of our center-manifold.

\subsection{Main result --- precise statement.}\label{s:main}

We are now in a position to state a precise version of Theorem \ref{thmCM} and Corollary \ref{corRE}. We are interested in  system \eqref{eqNLtime} and its modified variant,
\begin{equation}\label{e:nl}
\T u + \F(u)=0,
\end{equation}
\begin{equation}\label{e:nlm}
\T u + \F^\epsilon(u)=0.
\end{equation}
\begin{thm}[Center manifolds and reduced vector fields]\label{t:1}
Consider equations \eqref{e:nl} and \eqref{e:nlm} with assumptions (H1) on the linear convolution operator $\K$ and (H2) on the nonlinearity $\F$. Recall the definitions of the kernel $\E_0$ and the projection $\Q$ on $H^1_{-\eta}(\R,\R^n)$, \eqref{e:pro}. Then there exists a cut-off radius $\epsilon$, a weight $\delta>0$, and a map 
\[
\Psi:\ker\T\subset H^1_{-\delta}(\R,\R^n)\to \ker \Q\subset  H^1_{-\delta}(\R,\R^n),
\]
with graph 
\[
\M_0:=\left\{ u_0+\Psi(u_0) ~\middle|~ u_0 \in \ker \T \right\} \subset  H^1_{-\delta}(\R,\R^n),
\]
such that the following properties hold:
\begin{enumerate}
\item \emph{(smoothness)} $\Psi\in \cC^k$, with $k$ specified in (H2);
\item \emph{(tangency)} $\Psi(0)=0$, $D\Psi(0)=0$;
\item \emph{(global reduction)} $\M_0$ consists precisely of the solutions $u\in H^1_{-\delta}(\R,\R^n)$ of the modified equation \eqref{e:nlm};
\item \emph{(local reduction)} any solution $u\in H^1_{-\delta}(\R,\R^n)$ of the original equation \eqref{e:nl} with $\sup_{x\in\R}|u(x)|\leq \epsilon$ is contained in $\M_0$;
\item \emph{(translation invariance)} the shift $\tau_\xi,\xi\in\R$ acts on $\M_0$ and induces a flow $\Phi_\xi:\E_0\to\E_0$ through $\Phi_\xi=\Q\circ \tau_\xi\circ \Psi$;
\item \emph{(reduced vector field)} the reduced flow $\Phi_\xi(u_0)$ is of class $\cC^k$ in $u_0,\xi$ and generated by a reduced vector field $f$ of class $\cC^{k-1}$ on the finite-dimensional vector space $\E_0$. 
\end{enumerate}
In particular, small solutions on $t\in\R$ to $v'=f(v)$ on $\E_0$ are in one-to-one correspondence with small bounded solutions of \eqref{e:nl}.
\end{thm}

\paragraph{Higher regularity.} Completely analogous formulations of our main result are possible in spaces with higher regularity, $H^m_\eta(\R,\R^n)$, changing simply the assumptions on the nonlinearity, which will typically require higher regularity of pointwise nonlinearities, as  we shall see in Section \ref{s:non}. Moreover, one then concludes that small bounded solutions are in fact smooth in $x$, which one can, however, also conclude after using bootstrap arguments in the equation.

\subsection{Extensions --- parameters, symmetries, and pseudo-differential operators}\label{s:ext}

\paragraph{Parameters.}
In the context of bifurcation theory, one usually deals with parameter dependent problems. One then hopes to find center manifolds and reduced equations that depend smoothly on parameters. We therefore consider
\bqq
 u+\K \ast u +\F(u,\mu)=0,
\label{eqNLpar}
\eqq
for $u:\R\rightarrow\R^n$, $n\geq1$,  $\mu\in \R^d$, $d\geq 1$, and the nonlinear operator $\F$ is defined in a neighborhood of $(u,\mu)=(0,0)$. Again, we can define $\F^\epsilon=\F\circ (\chi_\epsilon,\mathrm{id})$, cutting off in the $u-$variable, only, leading to 
\bqq
 u+\K \ast u +\F^\epsilon(u,\mu)=0,
\label{eqNLparm}
\eqq

We then require a $\mu$-dependent version of Hypothesis (H2). 

\begin{Hypothesis}{(H2$_\mu$)}
We assume that there exists $k\geq 2$ and $\eta_0>0$ such that for all $\epsilon>0$, sufficiently small, the following properties hold.
\begin{enumerate}
\item $\F\in \cC^k(\V_u\times \V_\mu,W^{1,\infty}(\R,\R^n)$, for some small neighborhoods $\V_u\subset W^{1,\infty}(\R,\R^n)$, $\V_\mu\subset\R^d$, and $\F(0,0)=0$, $D_u\F(0,0)= 0$;
\item $\F$ commutes with translations for all $\mu$, $\F\circ\tau_\xi=\tau_\xi\circ\F$ for all $\xi\in\R$;
\item $\F^\epsilon:  H^1_{-\zeta}(\R,\R^n)\times \V_\mu \longrightarrow H^{1}_{-\eta }(\R,\R^n)$ is $\cC^k$ for all nonnegative pairs $(\zeta,\eta)$ such that $0<k\zeta <\eta<\eta_0$,  $D^j\F^\epsilon(u,\mu): (H^1_{-\zeta}(\R,\R^n))^j \longrightarrow H^{1}_{-\eta }(\R,\R^n)$ is bounded for $0<j\zeta\leq \eta<\eta_0$, $0\leq j\leq k$ and Lipschitz in $u$ for $1\leq j \leq k-1$, uniformly in $\mu \in \mathcal{V}_\mu$.
\end{enumerate}
\end{Hypothesis}

The analogue of the center manifold Theorem \ref{thmCM} for the parameter-dependent nonlocal equation \eqref{eqNLpar} is the following result.

\begin{thm}[Parameter-Dependent Center Manifold]\label{thmCMpar} 
Consider equations \eqref{eqNLpar}and \eqref{eqNLparm} with assumptions (H1) on the linear convolution operator $\K$ and  with assumption (H2$_\mu$) on the nonlinearity $\F$.  Recall the definition of kernel $\E_0$ and projection $\Q$ on $H^1_{-\eta}(\R,\R^n)$, \eqref{e:pro}. Then, possibly shrinking the neighborhood $\V_\mu$,  there exist a cut-off radius $\epsilon$, a weight $\delta>0$, and a map 
\[
\Psi:\ker\T\times \V_\mu\subset H^1_{-\delta}(\R,\R^n)\times\R^d\to \ker \Q\subset  H^1_{-\delta}(\R,\R^n),
\]
with graph 
\[
\M_0:=\left\{ (u_0+\Psi(u_0,\mu),\mu) ~\middle|~ u_0 \in \ker \T,\ \mu\in \V_\mu \right\} \subset  H^1_{-\delta}(\R,\R^n),
\]
such that the following properties hold:
\begin{enumerate}
\item \emph{(smoothness)} $\Psi\in \cC^k$, with $k$ specified in (H2)$_\mu$;
\item \emph{(tangency)} $\Psi(0,0)=0$, $D_{u_0}\Psi(0,0)=0$;
\item \emph{(global reduction)} $\M_0$ consists precisely of the pairs $(u,\mu)$, such that $u\in H^1_{-\delta}(\R,\R^n)$ is a solution of the modified equation \eqref{e:nlm} for this value of $\mu$;
\item \emph{(local reduction)} any pair $(u,\mu)$ such that $u$ is a solution $u\in H^1_{-\delta}(\R,\R^n)$ of the original equation \eqref{e:nl} with $\sup_{x\in\R}|u(x)|\leq \epsilon$ for this value of $\mu$ is contained in $\M_0$;
\item \emph{(translation invariance)} the shift $\tau_\xi,\xi\in\R$ acts on the $u$-component of $\M_0$ and induces a $\mu$-dependent flow $\Phi_\xi:\E_0\to \E_0$ through $\Phi_\xi=\Q\circ \tau_\xi\circ \Psi$;
\item \emph{(reduced vector field)} the reduced flow $\Phi_\xi(u_0;\mu)$ is of class $\cC^k$ in $u_0,\xi,\mu$ and generated by a reduced parameter-dependent vector field $f$ of class $\cC^{k-1}$ on the finite-dimensional vector space $\E_0$. 
\end{enumerate}
In particular, small solutions on $t\in\R$ to $v'=f(v;\mu)$ on $\E_0$ are in one-to-one correspondence with small bounded solutions of \eqref{e:nl}.
\end{thm}

\paragraph{Symmetries and reversibility}

In this subsection, we discuss the cases of equations possessing symmetries in addition to translation invariance. The aim is to show that such symmetries are inherited by the reduced equation. Generally speaking, we have an action of the direct product $G=\mathbf{O}(n)\times (\R\times\Z_2)$ on spaces of functions over the real line with values in $\R^n$, where $\mathbf{O}(n)$ is the group of orthogonal $n\times n$-matrices, and the action is defined through
\[
((\rho,\tau_\xi,\kappa)\cdot u)(x)=\rho\cdot u(\kappa(x-\xi)).
\]
Here, $\kappa x=-x$ when $\kappa$ is the nontrivial element of $\Z_2$. Note that $\chi_\epsilon$ commutes with the action of the full group  $\mathbf{O}(n)\times (\R\times\Z_2)$. 

\begin{Hypothesis}{(S)}
There is a subgroup $\Gamma\subset G$ that contains the pure translations, $\mathrm{id}\times\R\times\mathrm{id}\subset \Gamma$, such that \eqref{eqNLtime} is invariant under $\Gamma$, that is, 
\[
\gamma\circ \T=\T\gamma,\qquad \gamma\circ \F=\F\gamma,\qquad \mbox{ for all } \gamma\in\Gamma
.\]
We say the equation is reversible if $\Gamma\not\subset \mathbf{O}(n)\times\R\times\mathrm{id}$, that is, if the group of symmetries contains a reflection. We call $\Gamma_\mathrm{e}:=\Gamma\cap(\mathbf{O}(n)\times\R\times\mathrm{id})$ the equivariant part and $\Gamma_\mathrm{r}:=\Gamma\setminus \Gamma_\mathrm{e}$ the reversible part of the symmetries $\Gamma$.

\end{Hypothesis}
We remark that the equivariance properties of $\Q$ are concerned with symmetries in $\mathbf{O}(n)\times (\{0\}\times\Z_2)$, since the action of the shift on the kernel is induced through the projection itself, hence automatically respects the symmetry. 
We obtain the following result.

\begin{thm}[Equivariant Center Manifold]\label{thmCMequiv}
Assume that the above Hypotheses (H1), (H2), and (S) are satisfied. Then reduced center manifold $\M_0=\mathrm{graph}
\,(\Psi)$ and vector field $f$ from Theorem \ref{t:1} respect the symmetry, that is, 
\begin{enumerate}
\item $\E_0$ is invariant under $\Gamma$ and $\Q$ can be chosen to commute with all $\gamma\in\Gamma$;
\item $\Psi$ commutes with the action of $\Gamma$, $\M_0$ is invariant under the action of $\Gamma$;
\item $f$ commutes with the equivariant part,  $f\circ\gamma_1=\gamma_1\circ f$ for $\gamma=(\gamma_1,\tau_\xi,\mathrm{id})\in \Gamma_\mathrm{e}$,  and anti-commutes with the reversible part of the symmetries,   $f\circ\gamma_1=-\gamma_1\circ f$ for $\gamma=(\gamma_1,\tau_\xi,\kappa)\in \Gamma_\mathrm{r}$.
\end{enumerate}
\end{thm}
Analogous results hold for the parameter-dependent equation \eqref{eqNLpar}.

\paragraph{Pseudo-differential operators.}

Beyond operators of the form $\mathrm{id}+\K*$, one could consider more general nonlocal pseudo-differential operators and equations of the form
\bqq
 \mathscr{P} u+\F(u)=0,
\label{eqNLpd}
\eqq
for $u:\R\rightarrow\R^n$, $n\geq1$, where $\mathscr{P} u$ is a pseudo-differential operator defined as follows. Let $\nu \mapsto p(\nu)$ be an analytic function in $\mathscr{S}_{\eta_0}=\{|\Re(\nu)|\leq \eta_0\}\subset  \C$, and define 
\begin{equation*}
\mathscr{P} u(x):=\frac{1}{2\pi}\int_\R \rme^{\rmi \ell x} p(\rmi \ell) \widehat{u}(\rmi \ell)\md \ell, \quad \forall x\in \R,
\end{equation*}
with suitable assumptions on convergence of the integral, say, sufficient localization of $\hat{u}$. A typical assumption on $p$ requires asymptotic growth with fixed order $\alpha>0$,  $|p(\nu)- \nu^\alpha|\to 0$ for $|\nu|\to\infty$, together with derivatives. Assuming that $(p(\nu)-M)^{-1}$ is uniformly bounded in $\Sigma_{\eta_0}$ for some $M\in\R$, we can then precondition the equation as 
\bqq
(M-\mathscr{P})^{-1}\left( \mathscr{P} u+\F(u)\right)=-u+M(M-\mathscr{P})^{-1}u+(M-\mathscr{P})^{-1}\F(u)=0,
\label{eqNLpd2}
\eqq
which is of the form \eqref{e:nl}, with kernel given through $\widehat{\K}(\nu)=\frac{M}{M-p(\nu)}$. Kernel smoothness therefore is determined by the value of $\alpha$. It is worth noticing that this perspective allows us to construct center manifolds for higher-order differential equations without writing the equation as a first-order equation. 

The common feature of all those examples is that the leading-order part in the linearization, from a regularity point of view, is invertible, and the nonlinearity is bounded on the domain of the leading-order part. In those cases, ``preconditioning'' with the resolvent of the leading-order part gives an equation of the type considered here. From this perspective, forward-backward delay equation present an interesting extension, where the principal part is of the form 
\[
(\T u)(x)=\sum_j A_ju(x-\xi_j) + (\K*u)(x),
\]
for matrices $A_j$ and $\xi_j\in\R$, and a convolution kernel $\K$ with assumptions as considered earlier. Note that we do not include a derivative, as would be common for traveling-wave equations in lattices. Such forward-backward functional equations arise naturally when studying traveling waves in  space-time discretizations of partial differential equations. Of course, in some cases (in particular, when the $\xi_j$ are linearly dependent over $\mathbb{Q}$), the equation reduces to an equation over a lattice. Our approach can be applied whenever the characteristic equation of the  principal symbol
\[
d_0(\nu)=\mathrm{det}\,\left(\sum_j A_j\rme^{\nu\xi_j}\right),
\]
is invertible with uniform bounds on a complex strip $\mathscr{S}_{\eta_0}$. 

\subsection{Applying the result --- nonlinearities}\label{s:non}
Our goal here is to provide examples of nonlinearities that satisfy Hypothesis (H2) and more generally provide some basic tools that may help verifying (H2) in specific examples. 
We start with pointwise nonlinearities, then discuss nonlocal operators, and conclude with more general composition of operators.

\paragraph{Pointwise nonlinearities.} We first consider classical superposition operators, defined by pointwise evaluation of the composition. Let $g  \in \mathscr{C}^{k+1}(\R^n)$ for some $k \geq 2$ and define the superposition operator $\F$ as
\begin{equation*}
\F(u)[x]=g(u(x)), \quad \forall x \in \R, \quad \forall u \in W^{1,\infty}(\R,\R^n).
\end{equation*}
The properties listed in Hypothesis (H2) then are precisely the properties established in \cite[Lemma 3 \& 5]{vander-vangils:87}, with the small caveat that spaces $\cC^0_{-\eta}(\R,\R^n)$ instead of $H^1_{-\eta}(\R,\R^n)$ we considered there. Adapting the arguments is not difficult; we outline the key steps in Appendix \ref{a:non} for the convenience of the reader. 

One can also push these arguments to spaces $H^m_{-\eta}(\R,\R^n)$, requiring $g\in \cC^{k+m}(\R^n)$, thus allowing us to construct center manifolds in spaces $H^m_{-\eta}(\R,\R^n)$; see the remark after Theorem \ref{t:1}. We omit the details of this straightforward adaptation. 

\paragraph{Convolution operators.} A class of \emph{linear}  operators  that satisfies (H2), of course with the exception of $D\F(0)=0$, are convolutions with convolution kernel as in (H1). One can in fact generalize slightly, and consider convolutions $\K * u$ with kernel $\K$, an exponentially localized Borel measure, $\K=\K_0+\sum_j a_j \delta_{\xi_j}$, with $\K_0\in L^1_{\eta_0}(\R,\R^n)$ and $\sum_j|a_j|\rme^{\eta_0|\xi_j|}<\infty$.

\paragraph{Composition.} Slightly generalizing Hypothesis (H2), we can consider maps $\F$ mapping $\R^n$-valued functions to $\R^q$-valued functions, keeping all other properties from (H2). We claim that the composition of two such functions satisfying (H2) then also satisfies (H2), possibly with a smaller $\eta_0$. This can be readily obtained as follows. Consider the composition $\F\circ \G$ with derivative $\F'(\G(u))\cdot \G'(u)\cdot v$. We need to show that 
\[
\|\F(\G(u+v))-\F(\G(u))-\F'(\G(u))\cdot \G'(u)\cdot v\|_{H^1_{-\eta-\delta}(\R,\R^n)}=\co(\|v\|_{H^1_{-\eta}(\R,\R^n)}).
\]
For this, decompose as in the proof of the chain rule,
\begin{align*}
\F(\G(u+v))-\F(\G(u))-\F'(\G(u))\cdot \G'(u)\cdot v&=\left\{(\F(\G(u+v))-\F(\G(u))-\F'(\G(u))\cdot(\G(u+v)-\G(u))\right\}\\
&\quad+\left\{\F'(\G(u))\cdot(\G(u+v)-\G(u)-\G'(u)\cdot v)\right\}\\
&=:I+II.
\end{align*}
Now,
\[
\|I\|_{H^1_{-\eta-\delta}}=\co(\|\G(u+v)-\G(u)\|_{H^1_{-\eta}(\R,\R^n)})=\co(\|v\|_{H^1_{-\eta}(\R,\R^n)}),
\]
by differentiability of $\F$, and  Lipschitz continuity of $\G$. Next, 
\[
\|II\|_{H^1_{-\eta-\delta}(\R,\R^n)}=\co(\|v\|_{H^1_{-\eta}(\R,\R^n)}),
\]
by differentiability of $\G$ and boundedness of $\F'$. Boundedness of derivatives of the composition and Lipschitz continuity are readily checked form the chain rule formula. Higher derivatives are obtained in an analogous fashion; see also Appendix \ref{a:non} for similar arguments.

As a consequence, we can treat nonlinearities of the form $\K_1*f(\K_2*u,\K3*u,\ldots,\K_\ell*u)$, say. 

\paragraph{Multilinear convolutions.}
Slightly more general are multilinear convolution operators of the form 
\[
\mathscr{G}\cdot[u_1,\ldots,u_\ell](x)=\int\ldots\int \G(x_1-y_1,\ldots,x_\ell-y_\ell)u(y_1,\ldots,u(y_\ell)\md y_1\ldots \md y_\ell,
\]
with each components of $\G$ being in  $W^{1,1}_{\eta_0}(\R^\ell)$. One readily verifies that $\mathscr{G}$ is bounded as a multilinear operator on $H^1_{-\eta}(\R,\R^n)$. Again, convolution kernels generally  in the nonlinearity need not be smooth and may contain Dirac deltas.

%%%%%%%%%%%%%%%%%%%%%%%%%%%%%%%%%%%%%%%%%%%%%%%%%%%%%%%%%%%%%%%%%%%%%%%%%%%%%%%%%%%%%%%%%%%%%%%%%%%%%%%%%

\subsection{Applying the result --- projections}\label{s:pro}
The projection on the kernel clearly plays an important role in the actual computation of the reduced vector field $f$. We emphasize again that this projection cannot be canonically chosen as a spectral projection, as it acts on trajectories rather than a phase space. Abstractly speaking, projections on finite-dimensional subspaces always exist by Hahn-Banach's theorem.  In that respect, particular choices of projections could be favored over others mostly because they simplify computations: first, one would like to simplify the computation of the reduced vector field, and second, one would like to find good coordinates in which to analyze the reduced vector field. 

Since we are working in a Hilbert space $H^1_{-\delta}(\R,\R^n)$, one can of course simply use orthogonal projections on the kernel. In fact, since as we shall see later the kernel consists of smooth functions with at most polynomial growth, one can use a variety of weighted scalar products, and we shall briefly explore some choices below. On the other hand, we found it convenient in practical applications to use pointwise evaluations of functions and their derivatives as the arguably most easily computable projection.

To start with, we recall that a projection $\Q$ on a kernel $\mathrm{ker}\,(\mathcal{T})=\mathrm{span}\,(e_1,\ldots,e_M)$ can be identified with a collection of functionals $f^*_1,\ldots,f_M^*$ such that the Gram matrix $A$ with entries $A_{kl}=\langle e_l,f_k^*\rangle$ is invertible, by setting
\begin{equation}\label{e:gram}
Qu:=\sum_{j=1}^M \langle u,f_j^*\rangle e_j, \qquad \Q u:=A^{-1}Qu.
\end{equation}
In order to give specific examples, we need the following characterization of the kernel of $\mathcal{T}$. 
Recall the definition of the linear operator and its Fourier transform $\widehat{\T}(\nu):=I_n+\widehat{\K}(\nu)$, a matrix pencil defined and holomorphic on $\mathscr{S}_{\eta_0}=\left\{ \nu\in\C ~|~ |\Re(\nu)|<\eta_0 \right\}$. 
As a consequence, $d(\nu)=\det(\widehat{\T}(\nu))$ has finitely many roots, counted with multiplicity on the imaginary axis. Possibly reducing $\eta_0$, we assume that $d$ does not vanish off the imaginary axis and refer to roots as \emph{characteristic values}. We label those characteristic  values $\nu_j=\rmi \ell_j$, $1\leq j \leq m$, and denote by $r_j$ the dimension of $\ker \widehat{\T}(\nu_j)$, referred to in the sequel as  {\it geometric multiplicity} of the characteristic value $\nu_j$. Now let $e_{j,k}^0\in \C^n$, $1\leq k \leq r_j$, be a basis of the kernel,
\begin{equation}\label{e:chker}
\widehat{\T}(\nu_j)e_{j,k}^0=0. 
\end{equation}
Then there exist $n_{j,k} \geq r_j$ such that we can construct a maximal chain of root vectors $\left( e_{j,k}^p \right)_{0\leq p \leq n_{j,k}-1}$ which satisfy
\begin{equation}\label{e:chain}
\sum_{q=0}^p \dbinom{p}{q}\widehat{\T}^{(q)}(\nu_j)e_{j,k}^{p-q}=0, \quad 0\leq p \leq n_{j,k}-1, \quad \widehat{\T}^{(q)}(\nu_j)u:=\frac{d^q}{d\nu^q}\left(\widehat{\T}(\nu)u \right)_{|\nu=\nu_j},
\end{equation}
where $e_{j,k}^0:=e_{j,k}$. The sum $\alpha_j=n_{j,1}+\cdots+n_{j,r_j}$ is called the {\it algebraic multiplicity} of the characteristic value $\nu_j$. 

\begin{lem}\label{l:ker}
The maximal chain of root vectors is always finite and the algebraic multiplicity $\alpha_j$ coincides with the order of the root $\nu_j$ of $d(\nu)$.  Let $M\geq 1$ be defined as $M:=\alpha_1+\cdots+\alpha_m$. Then the kernel of $\T$ is isomorphic to $\R^M$, given explicitly through
\begin{equation*}
\E_0=\ker \T = \bigoplus_{j=1}^m \left( \bigoplus_{k=1}^{r_j} \mathrm{Span}\left\{ \vp_{j,k,p}(x), \ 0\leq p \leq n_{j,k}-1 \right\} \right), \qquad \vp_{j,k,p}(x) =\left(\sum_{q=0}^p \dbinom{p}{q}x^q e_{j,k}^{p-q}\right)\rme^{\rmi \ell_jx}.
\end{equation*}
\end{lem}

\begin{Proof}
The existence of Jordan chains as listed in \eqref{e:chain} is a standard result for analytic matrix pencils and can be proved readily using Lyapunov-Schmidt reduction on the eigenvalue problem; see for instance \cite{baum,gohberg} or 
\cite[Lem. 3.3]{jsw}. Now, taking the Fourier transform (in the sense of distributions) of $\T \varphi_{j,k,p}(x)=0$,  we readily find \eqref{e:chain}, thus showing that $u_0(x)$ indeed belongs to the kernel. Comparing dimensions, we find that the sums of the lengths of Jordan chains equals the multiplicity of the root of the determinant, again by standard theory for matrix pencils, we conclude that the elements $\varphi_{j,k,p}$ indeed form a basis of the kernel.
\end{Proof}

This particular basis gives us a representation of elements in the kernel in the form 
\begin{equation}\label{e:kerch}
u_0(x)= \sum_{j=1}^m \left(\sum_{k=1}^{r_j}\left(\sum_{p=0}^{n_{j,k}}A_{j,k,p} \left(\sum_{q=0}^p \dbinom{p}{q}x^qe_{j,k}^{p-q}\right)\rme^{\rmi \ell_jx}\right)\right)\in \E_0,
\end{equation}
and hence a canonical map  $\iota:\E_0\rightarrow \R^M$ through
\begin{equation*}
\iota(u_0)= \left\{ A_{j,k,p} \text{ for } 1\leq j \leq m, \quad 1\leq k \leq r_j , \quad 0 \leq p \leq n_{j,k} \right\}.
\end{equation*}
It is not difficult (but rather cumbersome in high-dimensional examples) to construct projections. We outline two possible choices. First, let $\omega(x)$ be a suitable weight function and define 
\[
(Qu)(x):=\left(\int_\R u(y)\overline{\vp_{j,k,p}(y)} \omega(y) \md y \right)\vp_{j,k,p}(x), \qquad \Q=A^{-1}Q,
\]
with $A$ as in \eqref{e:gram}. The entries of the Gram matrix $A_{kl}$ reduce to  integrals of the form 
$
\int_\R x^q\rme^{\rmi \ell x}\omega(x) \md x
$
which are explicitly given through derivatives of Gaussians and hyperbolic secants when $\omega(x)=\rme^{-x^2}$ or $\omega(x)=\mathrm{sech}(x)$, respectively. 
Note that, these projections, as the $L^2_{-\eta}(\R,\R^n)$-orthogonal projections, naturally extend to $H^m_{-\eta'}(\R,\R^n)$ for all $\eta'>0$, small enough, $m\geq0$. 

Second, in a different spirit, notice that the matrix $B=\left(b_{m,p}\right)_{1\leq m,p\leq n}$ with
\[
b_{m,p}:=\left\langle e^0,\left.\left(\partial_x-\rmi\ell\right)^{(m)}\right|_{x=0}\left(\sum_{q=0}^p \dbinom{p}{q}x^qe^{p-q}\rme^{\rmi\ell x}\right)\right\rangle=
\left\{\begin{array}{ll}
m! |e^0|^2,& m=p\\
0,& m>p,
\end{array}\right.
\]
is lower triangular with positive diagonal entries, hence invertible, thus yielding a canonical projection in the case of a simple root $\rmi\ell_j$ on sufficiently smooth functions $u\in H^N(\R,\R^n)$, $N$ sufficiently large. Generalizing to multiple roots is tedious but straightforward, taking additional derivatives when basis vectors $e_{j,k}^0$ and $e_{j',k'}^0$ are linearly dependent. While these projections are not defined on $H^1_{-\eta}(\R,\R^n)$ or even $L^2_{-\eta}(\R,\R^n)$, they can be used in computations whenever solutions are in fact smooth, typically because the nonlinearity maps into $H^m(\R,\R^n)$, locally. 

We note that the projections constructed here are defined for spaces $H^1_{-\eta}(\R,\R^n)$, say, for all weights $\eta>0$. Moreover, they commute with the natural embedding between those spaces. 

%%%%%%%%%%%%%%%%%%%%%%%%%%%%%%%%%%%%%%%%%%%%%%%%%%%%%%%%%%%%%%%%%%%%%%%%%%%%%%%%%%%%%%%%%%%%%%%%%%%%%%%%%

\subsection{Applying the result --- Taylor jets}\label{subsec:example}

We will apply our main result later on but want to give a fairly trivial example of how to compute Taylor jets in practice, here. In fact, the  procedure of deriving the reduced system \eqref{eqReduced} involves algebra that is somewhat different from the more commonly known algebra associated with Taylor jets in phase space and ordinary center manifolds. We consider a scalar nonlocal equation of the form,
\bqq
\label{ex1}
u+\K*u -u^2 =0,
\eqq
where we suppose that $\K$ satisfies Hypothesis (H1) for a given $\eta_0>0$ together with the assumptions that
\begin{equation*}
\quad \int_\R \K(x) \md x =-1, \quad \int_R x\K(x)\md x=-\alpha^{-1} \neq 0, \quad \text{ and } \quad d(\rmi\ell)=1+\widehat{\K}(\rmi\ell)\neq0 \text{ for all } \ell \in\R\setminus\{0\}.
\end{equation*} 
As a consequence,   $\E_0=\ker \T= \left\{ 1 \right\}$,  the constant functions. A natural candidate for the projection onto the kernel is $(\Q u)(x) \equiv  u(0) \in \E_0$, clearly defining a bounded  projection on $H^1_{-\eta}$ onto $\E_0$ for any $0<\eta<\eta_0$. Furthermore, the nonlinear operator $\F(u)=-u^2$ is a Nemytskii operator and satisfies Hypothesis (H2) as discussed above. Our main result, Theorem \ref{thmCM}, then implies existence of a  center manifold $\M_0$, and any small bounded solutions of \eqref{ex1} can be written as
\begin{equation*}
u = u_0 + \Psi(u_0),
\end{equation*}
where $u_0:=A\cdot 1 \in \E_0$. As the map $\Psi$ is $\cC^k$ for any $k\geq 2$, we can look for its Taylor expansion near $0$, and using the properties $\Psi(0)=D_u\Psi(0)=0$, we obtain
\begin{equation*}
\Psi(u_0) = A^2 u_1+A^3 u_2+ \mathcal{O}(A^4).
\end{equation*}
Inserting this ansatz into the nonlocal equation \eqref{ex1} and identifying terms of order $A^2$, we obtain that $u_1$ should satisfy
\begin{equation*}
\T u_1 = 1, \text{ with } \Q(u_1)=0. 
\end{equation*}
 Using that $\int x\K(x)\md x\neq 0$, we 
obtain that  $u_1(x) = \alpha x$, for all $x\in \R$. At cubic order,  we find that
 \begin{equation*}
 \T u_2 = 2 u_1, \text{ with } \Q(u_2)=0. 
 \end{equation*}
 We look for solution $u_2$ that can be written as $u_2(x) = \beta_2 x^2+\beta_1 x$, which leads to the compatibility conditions
 \begin{align*}
 \beta_2 \int_\R \K(y) y^2 \md y +\frac{\beta_1}{\alpha} &=0,\\
 2\frac{\beta_2}{\alpha}&=2\alpha, 
 \end{align*}
 such that $\beta_2=\alpha^2$ and $\beta_1:=- \kappa_2 \alpha ^3$, where $\kappa_2:=\int_\R \K(y) y^2 \md y$. Finally, we apply the definition of the flow to our solution
\begin{equation*}
u(x) = A + \alpha xA^2 +(\alpha^2x^2- \kappa_2\alpha^3x)A^3+ \mathcal{O}_x(A^4),
\end{equation*}
to obtain that
\begin{align*}
\varphi_x(A) &= \Q \left[ A + \alpha (\cdot + x)A^2 + \left(\alpha^2(\cdot + x)^2-\kappa_2\alpha^3(\cdot + x)\right)A^3 +\mathcal{O}_{(\cdot + x)}(A^4) \right] \\
&= A+\alpha x A^2  +  \left(\alpha^2x^2-\kappa_2\alpha^3x\right)A^3 +\mathcal{O}_x(A^4).
\end{align*}
Given a smooth flow, we obtain a vector field in the standard fashion, differentiating the flow at time $x=0$,
\begin{equation*}
\frac{\md \varphi_x}{\md x}|_{x=0}= \alpha A^2-\kappa_2\alpha^3A^3  + \mathcal{O}(A^4),
\end{equation*}
thus giving the Taylor expansion of the reduced equation up to third order through
\begin{equation*}
\frac{\md A}{\md x}= \alpha A^2-\kappa_2\alpha^3A^3   + \mathcal{O}(A^4).
\end{equation*} 
This strategy of differentiating the flow induced by the shift of bounded solutions, pulled back to the kernel, is at the heart of our construction of reduced vector fields in the proof of Theorem \ref{t:1} in Section \ref{subsec:smooth}.

Note that, absent further parameters,  the reduced differential equation, here,  does not possess any non-trivial bounded solutions. In other words, the center manifold here yields a uniqueness result for small bounded solutions, in a class of sufficiently smooth functions. Adding parameters, one would find the typical heteroclinic trajectories in a saddle-node bifurcation.

%%%%%%%%%%%%%%%%%%%%%%%%%%%%%%%%%%%%%%%%%%%%%%%%%%%%%%%%%%%%%%%%%%%%%%%%%%%%%%%%%%%%%%%%%%%%%%%%%%%%%%%%%

\section{Proofs of the main results}\label{s:proofs}
In this section, we give proofs of our main results. We start with the characterization of the kernel and the analysis of the linearization in exponentially weighted spaces in Section \ref{s:pker}. We prove existence and regularity of the center manifold in Section \ref{s:pf}, following very much the standard approach via contraction mapping principles on scales of Banach spaces. Section \ref{subsec:smooth} establishes smoothness of the flow on the kernel induced by translations of bounded solutions via bootstraps and thereby establishes existence of a reduced vector field governing the set of bounded solutions. 
Finally, Section \ref{s:sym} outlines modifications and adaptations in the cases with additional symmetries. 

%%%%%%%%%%%%%%%%%%%%%%%%%%%%%%%%%%%%%%%%%%%%%%%%%%%%%%%%%%%%%%%%%%%%%%%%%%%%%%%%%%%%%%%%%%%%%%%%%%%%%%%%%

\subsection{Properties of the linearization}\label{s:pker}

We give characterizations of bounded solutions of the linear part of our equation and establish bounded invertibility of a suitably bordered equation. 

Consider therefore the linearization 
\begin{equation}\label{e:lin}
{\T}:H^1_{-\eta}(\R,\R^n) \longrightarrow H^{1}_{-\eta}(\R,\R^n),\qquad \T u=-u+\K*u,\qquad 0<\eta\ll 1,
\end{equation}
with associated characteristic equation $d(\nu):=\mathrm{det}(\hat{\T}(\nu))$. 

\begin{lem}\label{l:linf}
The operator $\T$ defined in \eqref{e:lin} is Fredholm of index $M$ and onto, where $M$ is the sum of the multiplicities of roots of $d(\nu)$ on $\nu\in\rmi\R$. 
\end{lem}
\begin{Proof}
This result is a direct consequence of \cite{faye-scheel:13}, in particular Theorem 3 and Lemma 5.1 from this reference.
Since in this reference, we considered matrix operators of the form $\frac{\md }{\md x}+A+K*$, we first convert our operator into this form, writing  $\T= \cD^{-1} \left(\cD\T \right)$, where 
$\cD: H^1_{-\eta}(\R,\R^n)\longrightarrow L^2_{-\eta}(\R,\R^n)$, $u\mapsto \frac{\md u}{\md x}+\rho u$ is an isomorphism provided $\rho>\eta_0$, thus reducing the problem to establishing Fredholm properties of $\cD \T$, which is of the form $\cD \T:= \frac{\md }{\md x}+\cN$, where $\cN(u):=\left(\K'+\rho\K +\rho\delta_0 \right)*u$ and $\delta_0$ is the Dirac delta function. 

Fredholm properties of operators  such as $\cD \T$ have been studied in \cite{faye-scheel:13} where it was shown that $\cD \T: H^1_{-\eta}(\R,\R^n)\longrightarrow L^2_{-\eta}(\R,\R^n)$ is a Fredholm operator \cite[Theorem 2]{faye-scheel:13} with index $\dim \E_0$  \cite[Corollary 4.9]{faye-scheel:13}. Roughly speaking, one conjugates the operator with the multiplier $\cosh(\eta x)$ to find an $x$-dependent convolution operator of the form considered in this reference. 

Theorem 3 and Corollary 4.9 of \cite{faye-scheel:13} state that the Fredholm index is given by the spectral flow, in this case, the number of roots of the characteristic equation on the imaginary axis, counted with multiplicity. Since  the characteristic equation associated to $\cD\T$ is given by
\begin{equation*}
(\nu+\rho)^n\det\left(I_n+\widehat{\K}(\nu)\right)=0,
\end{equation*}
roots on the imaginary axis stem from roots of $d(\nu)$, only, which proofs the result.
\end{Proof}
We now augment equation \eqref{eqNL} with the ``initial condition'', $\Q(u)=u_0$, for a given parameter  $u_0\in \E_0$, which leads us to consider the ``bordered'' operator
\bqq
\label{eqker}
\begin{array}{rrcl}
\widetilde{\T}: &H^1_{-\eta}(\R,\R^n) &\longrightarrow& H^{1}_{-\eta}(\R,\R^n)\times\E_0\\
& u & \longmapsto & \left( \T(u), \Q(u)\right).
\end{array}
\eqq 

\begin{lem}\label{l:21}
For any $0<\eta<\eta_0$, $\widetilde{\T}$ defined in \eqref{eqker} is invertible with bounded inverse,
\bqq
\| \widetilde{\T}^{-1}\|_{H^{1}_{-\eta }(\R,\R^n) \rightarrow H^1_{-\eta}(\R,\R^n)\times \E_0} \leq C(\eta),
\label{invT}
\eqq
with $C(\eta)<\infty$ continuous for $0<\eta<\eta_0$.
\end{lem}

\begin{Proof}
Since we are adding finitely many dimensions to the range, Fredholm bordering implies that $\widetilde{\T}$ is Fredholm, of index 0. Whenever $\widetilde{\T}u=0$, we conclude that $\T u=0$ from the first component, hence $u\in \E_0$. The second component implies that $\Q(u)=0$, which for $u\in \E_0$ implies $u=0$. 
\end{Proof}

%%%%%%%%%%%%%%%%%%%%%%%%%%%%%%%%%%%%%%%%%%%%%%%%%%%%%%%%%%%%%%%%%%%%%%%%%%%%%%%%%%%%%%%%%%%%%%%%%%%%%%%%%

\subsection{Lipshitz and smooth center manifolds}\label{s:pf}
We now rewrite equations \eqref{eqNL} together with \eqref{eqker}, using the modified nonlinearity $\F^\epsilon$ instead of $\F$,  into a more compact form
\bqq
\label{eqNLnew}
\widetilde{\T}(u)+\widetilde{\F^\epsilon}(u;u_0)=0,
\eqq 
where 
\[
\widetilde{\F^\epsilon}(u;u_0)=(\F^\epsilon(u),-u_0).
\]
Applying $\widetilde{\T}^{-1}$ to equation \eqref{eqNLnew}, we obtain an equation of the form
\bqq
\label{eqNLeps}
u=-\widetilde{\T}^{-1}\left(\widetilde{\F}^\epsilon(u;u_0)\right):=\cS^\epsilon(u;u_0),
\eqq
for any $u_0\in \E_0$.  We view \eqref{eqNLeps} as a fixed point equation with parameter $u_0$ and establish that  that $\cS^\epsilon(\cdot;u_0)$ is a contraction map on $H^1_{-\eta}(\R,\R^n)$. From the definition of $\F^\epsilon$ and the fact that $\F^\epsilon(0)=D\F^\epsilon(0)=0$ with $\F^\epsilon$ of class $\cC^k$ for $k\geq 2$ on $W^{1,\infty}(\R,\R^n)$, one obtains the following estimates as $\epsilon\rightarrow0$, 
\begin{subequations}
\begin{align}
\delta_0(\epsilon) &:= \underset{u \in H^1_{-\eta}(\R,\R^n)}{\sup}\| \F^\epsilon(u) \|_{H^{1}_{-\eta }(\R,\R^n)}= \mathcal{O}(\epsilon^2), \label{estF}\\
\delta_1(\epsilon)&:=\text{Lip}_{H^{1}_{-\eta }(\R,\R^n)}(\F^\epsilon)=\mathcal{O}(\epsilon). \label{estLip}
\end{align}
\end{subequations}
Indeed, by definition, we have $\F^\epsilon(u)(x)=\F(u)(x)$ whenever $\| u(x) \| \leq \epsilon $ and $\F^\epsilon(u)(x)=0$ whenever $\| u(x) \| \geq 2\epsilon$. Using the fact that $H^1$ functions are also continuous functions, we obtain the desired estimates by further noticing that $\F^\epsilon(u)$ is superlinear near $u=0$.
In turn, these estimates imply
\begin{align*}
\| \cS^\epsilon (u;u_0)\|_{H^1_{-\eta}(\R,\R^n)} &\leq C(\eta)\left(\delta_0(\epsilon)+ \|u_0\|_{H^1_{-\eta}(\R,\R^n)}\right),\\
\| \cS^\epsilon (u;u_0)-\cS^\epsilon (v;u_0) \|_{H^1_{-\eta}(\R,\R^n)} &\leq C(\eta) \delta_1(\epsilon) \| u-v \|_{H^1_{-\eta}(\R,\R^n)}, 
\end{align*}
for all $u,v\in {H^1_{-\eta}(\R,\R^n)}$ and $u_0\in \E_0$. Let $\bar \eta \in (0,\eta_0)$ and $\tilde \eta \in (0, \bar \eta /k)$, then, for sufficiently small $\epsilon$, we have 
\begin{equation*} C(\eta) \delta_1(\epsilon)<1,\quad \forall \eta \in [\tilde \eta, \bar \eta].
\end{equation*} 
As a consequence, there exists a unique fixed point $u=\Phi(u_0) \in {H^1_{-\eta}}(\R,\R^n)$. From Lipshitz continuity of the fixed point iteration, we conclude that  $\Phi$ is a Lipschitz map, and $\Phi(0)=0$ by uniqueness of the fixed point. For each $\eta \in [\tilde \eta, \bar \eta]$, this defines a continuous map $\Psi:\E_0\rightarrow \ker \Q \subset H^1_{-\eta}(\R,\R^n)$ so that
\begin{equation*}
u=\Phi(u_0):=u_0+\Psi(u_0).
\end{equation*} 

\begin{lem}\label{l:cmfdsm}
Under the Hypotheses (H1)-(H3) we have for each $p$ with $1\leq p \leq k$ and for each $ \eta \in (p\tilde\eta,\bar\eta]$ that $\Psi:\E_0\rightarrow H^1_{-\eta}(\R,\R^n)$ is of class $\cC^p$.
\end{lem}

\begin{Proof} First, notice that $\Phi$ shares the same properties as $\Psi$ so that it is enough to prove the Lemma for the map $\Phi$. We also recall that the modified nonlinearity $\F^\epsilon$ is $\cC^k$ from $H^1_{-\zeta} (\R,\R^n)$ to $ H^{1}_{-\eta }(\R,\R^n)$ for any $\zeta$ and $\eta$ satisfying $0<k\zeta <\eta<\eta_0$. Furthermore, we have that $D^j\F^\epsilon(u): (H^1_{-\zeta}(\R,\R^n))^j \longrightarrow H^{1}_{-\eta }(\R,\R^n)$ is bounded for $0<j\zeta\leq \eta<\eta_0$, $0\leq j\leq k$ and Lipschitz in $u$ for $1\leq j\leq  k-1$. The regularity properties of $\F^\epsilon$ are automatically inherited by $\cS^\epsilon$ by boundedness of the map $\widetilde{\T}^{-1}$. The conclusion of the lemma is then an application of the contraction mapping theorem on scales of Banach spaces as presented in \cite{vander-vangils:87}. The adaptations are straightforward; the main steps are outlined in Appendix \ref{secapp:thm}.
% that we state in the Appendix \ref{secapp:thm} in Theorem \ref{thm:appendix} for convenience. We now explain the main steps.\\

\end{Proof}

%%%%%%%%%%%%%%%%%%%%%%%%%%%%%%%%%%%%%%%%%%%%%%%%%%%%%%%%%%%%%%%%%%%%%%%%%%%%%%%%%%%%%%%%%%%%%%%%%%%%%%%%%

\subsection{Smoothness of the reduced flow and reduced vector fields}\label{subsec:smooth}

In this subsection, we establish that the flow on the center manifold is smooth such that we can obtain the reduced ordinary differential equation \eqref{eqReduced} simply through differentiating the flow at time zero. Consider the action of the shift operator on functions, defined through
\bqq
\label{eqphi}
\begin{array}{lcl}
\R\times H^{1}_{-\eta}(\R,\R^n) & \longrightarrow & H^1_{-\eta}(\R,\R^n) \\
(x ,u) &\longmapsto &\phi(x,u) := u(\cdot+x),
\end{array}
\eqq
for any $0<\eta<\eta_0$. We briefly write $\phi_x:=\phi(x,\cdot):H^{1}_{-\eta}(\R,\R^n) \longrightarrow L^2_{-\eta}(\R,\R^n)$. Clearly, $\phi_x$ is  bounded linear. 
Therefore, and by translation invariance of the original equation, $\phi_x$ maps bounded solutions to bounded solutions. The following commutative diagram shows how this action of the shift induces a flow on the kernel $\E_0$, 
\begin{center}
$\xymatrix @!=2.5cm{
\E_0 \ar[r]^{\mathrm{id}+\Psi} \ar@{}[dr]|*{\fontsize{15cm}{15.5cm}\selectfont \circlearrowright} \ar[d]_*{\varphi_x} &  H^1_{-\eta}(\R,\R^n) \ar[d]^*{\phi_x} \\
\E_0 \ar@<3pt>[r]^{\mathrm{id}+\Psi} & H^1_{-\eta}(\R,\R^n) \ar@<3pt>[l]^{{\Q}}}
$  
\hfill
$\xymatrix @!=2.5cm{
\E_0 \ar[r]^{\mathrm{id}+\Psi} \ar@{}[dr]|*{\fontsize{15cm}{15.5cm}\selectfont \circlearrowright} \ar[d]_*{\varphi_x} &  H^1_{-\eta}(\R,\R^n) \ar[d]^*{\iota\circ\phi_x} \\
\E_0 \ar@<3pt>[r]^{\iota\circ(\mathrm{id}+\Psi)} & L^2_{-\eta}(\R,\R^n) \ar@<3pt>[l]^{\widetilde{\Q}}}
$
\hfill
$\xymatrix @!=2.5cm{
\E_0 \ar[r]^{\iota\circ(\mathrm{id}+\Psi)} \ar@{}[dr]|*{\fontsize{15cm}{15.5cm}\selectfont \circlearrowright} \ar[d]_*{\varphi_x} &  L^2_{-\eta}(\R,\R^n) \ar[d]^*{\iota\circ\phi_x\circ\iota^{-1}} \\
\E_0 \ar@<3pt>[r]^{\iota\circ(\mathrm{id}+\Psi)} & L^2_{-\eta}(\R,\R^n) \ar@<3pt>[l]^{\widetilde{\Q}}}
$
\end{center}

The left diagram, $\mathrm{id}+\Psi$ denotes the parameterization of bounded solutions over the kernel. On the right, $\phi_x$ denotes the shift which is pulled back to the kernel via the projection $\Q$, the inverse of $\mathrm{id}+\Psi$. The right diagram views the bounded solutions as elements of $L^2_{-\eta}(\R,\R^n)$, by composing the parameterization $\mathrm{id}+\Psi$ with the embedding $\iota:H^1_{-\eta}(\R,\R^n)\to L^2_{-\eta}(\R,\R^n)$. The inverse of the parameterization is the extension of the projection $\Q$ to $L^2_{-\eta}(\R,\R^n)$. The induced flow on the kernel $\E_0$ is naturally the same as in the left diagram. In the center diagram, we view the shift as a map from $H^1_{-\eta}(\R,\R^n)$ into $L^2_{-\eta}(\R,\R^n)$. Clearly, $\iota\circ{\Phi_x}$ is continuously differentiable in $x$, with derivative given by the bounded linear map $\frac{\md y}{\md x}$.  Since $\widetilde{\Q}$ is a bounded projection on $L^2_{-\eta}(\R,\R^n)$, we find that 
\[
\varphi_x:=\tilde{Q}\phi_x\circ(\mathrm{id}+\Psi),
\]
is continuously differentiable in $x$. 
From Theorem \ref{thmCM} we know that $\Psi$ is a $\cC^k$ map from $\E_0$ to $H^1_{-\eta}(\R,\R^n)$. Therefore, the map $x\mapsto \varphi_x$ inherits the regularity properties of $\phi$, from which we deduce that $\frac{\md \varphi_x}{\md x}|_{x=0}$ is a $\cC^k$ vector field on $\E_0$, 
\bqq
\label{defReduced}
\frac{\md \varphi_x}{\md x}|_{x=0}=:f(u_0).
\eqq
Conversely, solutions to $\frac{\md u}{\md x}=f(u)$, $u(0)=u_0$ yield trajectories $\varphi_x(u_0)$ and solutions to the nonlocal equation $(\mathrm{id}+\Psi)(\varphi_x(u_0))$.

\subsection{Proof of Theorem \ref{t:1}}

We conclude the proof of Theorem \ref{t:1}. We established in Section \ref{s:pf} the existence of the map $\Psi$ and the assosiated smooth manifold $\M_0$. By uniqueness, and since $\F(0)=0$ implies that $u(x)\equiv 0$ is a solution, $\Psi(0)=0$. Differentiating the equation \eqref{e:nlm} at $u=u_0+\Psi(u_0)$ with respect to $u_0$ at $u_0=0$ gives that $D\Psi(0)=0$ viewed as an operator from $H^1_{-\eta}(\R,\R^n)$ to $H^1_{-\eta-\delta}(\R,\R^n)$ for any $\delta>0$, which implies that the derivative as a map from $H^1_{-\eta}(\R,\R^n)$ into itself also vanishes, which establishes (ii). Global reduction (iii) is a consequence of the construction as a contraction mapping, ensuring a unique fixed point for any $u_0\in\E_0$. Translation invariance and the existence of a reduced vector field, properties (v) and (vi), were discussed in Section \ref{subsec:smooth}. Local reduction, property (iv), follows since the nonlinearity is identical to the modified nonlinearity on the ball of size $\varepsilon$. It remains to show that small solutions to the reduced differential equation yield solutions to the original problem. To see this, notice that smallness of the trajectory in $\E_0$ implies, by construction of the flow and continuity of the map $\Psi$, smallness of all translates of the solution $u(x)$ in $H^1_{-\delta}(\R,\R^n)$, which readily establishes smallness in $L^\infty$ and concludes the proof of Theorem \ref{t:1}.

\subsection{Symmetries and parameters --- proofs}\label{s:sym}

We conclude the proofs of our main results by addressing the extensions in Theorem \ref{thmCMpar} and  \ref{thmCMequiv}.

 \begin{Proof}[ of Theorem \ref{thmCMpar}.]
 We cast the parameter-dependent system \eqref{eqNLpar} as a particular case of \eqref{eqNL}, in the form,
 \bqq
\bu+\J \ast \bu +\cR(\bu)=0,
\label{eqNLbis}
\eqq
by setting $\bu:=(u,\mu)$, and 
 \begin{align*}
 \B&:=\left(\begin{matrix} I_n & D_\mu\F(0,0) \\ 0_{1,n} & 1 \end{matrix}\right),\\
 \J &:=\B^{-1}\left(\begin{matrix} \K & 0_{n,1} \\ 0_{1,n} & \I \end{matrix}\right),\\
 \cR(\bu)&:=\B^{-1}(\F(u,\mu)-D_\mu\F(0,0)\mu,0),
 \end{align*}
where  $\I:= -\left(1+\frac{\md}{\md x} -\frac{\md^2}{\md x^2}\right)^{-1}$. Indeed, we first use the fact that $\mu$ is a parameter such that
 \begin{equation*}
 -\mu+\mu+\frac{\md \mu}{\md x}-\frac{\md^2 \mu}{\md x^2} =0.
 \end{equation*}
Applying the convolution operator $(1+\frac{\md}{\md x}-\frac{\md^2}{\md x^2})^{-1}$, we obtain 
 \begin{equation*}
 \mu - \left(1+\frac{\md}{\md x}-\frac{\md^2}{\md x^2} \right)^{-1} \mu =0,
 \end{equation*}
 which can be cast as the nonlocal equation
 \begin{equation*}
\mu + \I*\mu =0.
 \end{equation*}
One readily finds that  $\I(x) \in W_\alpha^{1,1}$ for $|\alpha|<(\sqrt{5}-1)/2)$. As consequence,  Hypothesis (H1) is satisfied for $\J$. Furthermore, it is clear that Hypothesis (H2$_\mu$) for $\F$ in \eqref{eqNLpar} implies that Hypothesis (H2) is satisfied for $\cR$. Since all solutions necessarily have $\mu(x)$ constant in $x$, this proves the theorem. 
  \end{Proof}

\begin{Proof}[ of Theorem \ref{thmCMequiv}.]
First notice that the cut-off, performed with respect to the norm in $\R^n$ which is invariant under the action of the orthogonal group, preserves equivariance as stated in Hypothesis (S). 
The uniqueness of the fixed point of the equation \eqref{eqNLeps} in the proof of Theorem \eqref{thmCM} implies that the corresponding center manifold is invariant under $\mathbf{S}$, provided that equation \eqref{eqNLeps} is equivariant under $\mathbf{S}$. Since the convolution part of $\T$ is equivariant with respect to $\mathbf{S}$ so will be $\T$ and the projection $\Q$, and thus $\widetilde{\T}$ is also equivariant. The properties of $f$ follow from differentiation of the properties of the flow. 
\end{Proof}

\section{Applications}\label{sec:appli}

We describe two applications of our center-manifold result to questions of existence of coherent structures in neural field equations. We construct  stationary solutions and traveling waves as examples in Sections \ref{s:stat} and \ref{s:tw}, respectively.  The emphasis here is on illustrating the feasibility of the reduction and the mechanics of the computation rather than  motivation for the problems or techniques to analyze reduced equations.

\subsection{Stationary solutions of neural field equations -- mode interactions}\label{s:stat}

We study small bounded solutions of neural field equations \eqref{nfe}, which take the form:
\bqq
0= -u+\K \ast S(u,\mu),
\label{Nfe}
\eqq
for some bifurcation parameter $\mu >0$. Such problems have been investigated in the literature either numerically or for very specific kernels with rational Fourier transform; see \cite{faye-etal:13} and references therein. In the above equation, $u:\R \rightarrow \R$ is a scalar unknown, the kernel function $\K$ and the nonlinearity $S$ satisfy the hypotheses below, reflecting simple modeling assumptions. We refrain from exploring minimal regularity assumptions on the nonlinearity and work with  smooth functions. Also, to avoid overly involved computations, we restrict ourselves to odd nonlinearities, in particular precluding quadratic terms in the Taylor jet of the center manifold. We also restrict  to the most relevant class of symmetric kernels.

\begin{hyp}
We suppose that the nonlinear function $S$ satisfies the following properties:
\begin{itemize}
\item[(i)] $(u,\mu)\mapsto S(u,\mu)$ is smooth on $\R^2$ with $|S(u,\mu)| \leq s_m$ and $0 \leq D_u S(u,\mu) \leq \mu s_m$ for all $(u,\mu)\in \R \times (0,+\infty)$ for some $s_m>0$;
\item[(ii)] $u\mapsto S(u,\mu)$ is an odd function, and 
\begin{equation*}
S(u,\mu) = \mu u -\frac{u^3}{3}+\mathcal{O}(|u|^5), \text{ as } u \rightarrow 0,
\end{equation*}
for all $\mu>0$.
\end{itemize}
\end{hyp}

\begin{hyp}
Let $\eta_0>0$. We suppose that $\K\in W_{\eta_0}^{1,1}(\R)$ is symmetric. Furthermore, we assume that the characteristic equation $d(\nu,\mu) = -1+ \mu  \widehat{\K} (\nu)$ satisfies:
\begin{itemize}
\item[(i)] $d(\nu,\mu)=\left[-\left( \nu^2+\ell_c^2\right)^2+\mu-\mu_c \right] \widetilde{d}(\nu,\mu)$ for a unique  $(\ell_c,\mu_c)\in(0,+\infty)^2$ such that  $\mu_c \widehat{\K} (\rmi \ell_c)=1$;
\item[(ii)] the function $\nu \mapsto  \widetilde{d}(\nu,\mu)$ does not have any roots on the imaginary and is analytic in the strip $\mathscr{S}:=\left\{ \nu\in\C ~|~ |\Re(\nu)|<\eta_0 \right\}$ for all $\mu>0$.
\end{itemize} 
\end{hyp}

\paragraph{Notation.} For any $(m_1,m_2)\in \N \times \Z$, we denote
\bqq
\kappa_{m_1,m_2}:=\int_\R x^{m_1} \K(x) \rme^{-m_2 \rmi \ell_c x}\md x.
\label{eq:kappa}
\eqq
From our condition on the characteristic equation, we have that  $\kappa_{0,\pm1}=1/\mu_c$ and $\kappa_{1,\pm1}=0$. From the symmetry of the kernel $\K$, we have that if $m_1\in\N$ is even, then $\kappa_{m_1,m_2}=\kappa_{m_1,-m_2}=\overline{\kappa}_{m_1,m_2}\in \R$, and if $m_1\in\N$ is odd, then $\kappa_{m_1,m_2}=-\kappa_{m_1,-m_2}=-\bar\kappa_{m_1,m_2}\in \rmi \R$.

With these hypotheses in hand, we define two usual operators
\begin{align*}
\T u &:= -u + \mu_c \K \ast u,\\
\F(u,\lambda)&:= \K \ast \left[ S(u, \lambda+\mu_c) - \mu_c  u\right],
\end{align*}
such that equation \eqref{Nfe} can be written as
\bqq
0 = \T u + \F(u,\mu-\mu_c).
\label{NFE}
\eqq

\paragraph{Symmetries.} It is important to notice that, in addition to the translation equivariance, equation \eqref{NFE} possesses two other symmetries, that we denote by $\mS_1$ and $\mS_2$ respectively, acting on functions via
\begin{equation*}
\mS_1 u(x) := u(-x), \quad \text{and }\quad \mS_2 u(x) :=-u(x), \quad \forall x \in \R.
\end{equation*}
The first symmetry is a consequence of the fact that the kernel $\K$ is a symmetric function, whereas the second symmetry results from the odd symmetry of the nonlinearity $S$ with respect to its first argument. Finally, let us remark that the conditions on the dispersion relation ensures that the kernel $\E_0$ of $\T$ is given by
\begin{equation*}
\E_0 =\text{Span} \left\{ \rme^{\pm \rmi \ell_c x}, x \rme^{\pm \rmi \ell_c x} \right\} \subset H^1_{-\eta}(\R),
\end{equation*}
for all $0<\eta<\eta_0$. In the following, we shall denote $\zeta_0(x):=\rme^{ \rmi \ell_c x}$ and  $\zeta_1(x):=x\rme^{ \rmi \ell_c x}$, with $\bar \zeta_0$ and $\bar \zeta_1$ their respective complex conjugate. As a consequence, any functions $u_0\in \E_0$, can be decomposed as
\bqq
\label{decomposition}
u_0 = A \zeta_0+\overline{A \zeta_0}+B \zeta_1+\overline{B \zeta_1}\in \E_0,
\eqq
for $(A,B)\in \C^2$. We remark that the actions of $\mS_{1,2}$ on $u_0$ are given by 
\begin{align*}
\mS_1 u_0 &= \bA \zeta_0+A \overline{\zeta_0} - \B \zeta_1 -B  \overline{\zeta_1}, \\
\mS_2 u_0 &= -A \zeta_0-\overline{A \zeta_0}-B \zeta_1-\overline{B \zeta_1}.
\end{align*}
We identify the action of $\mS_{1,2}$ on the quadruplet $(A,\bA,B,\bB)$ as
\begin{align*}
\mS_1 \cdot(A,\bA,B,\bB) &= (\bA, A , - \bB , -B  ), \\
\mS_2 \cdot(A,\bA,B,\bB) &= (-A,-\bA,-B -\bB).
\end{align*}

\paragraph{Projection $\Q$.} We now define the projection $\Q$ from $H^4_{-\eta}\rightarrow \E_0$. Note that by Sobolev embedding we have $H^4(\R)\subset \mathscr{C}^3(\R)$, and thus we can take linear combinations of $u^k(0)$ for any $k=0,\ldots,3$. For any $u_0\in \E_0$ that can be written as in \eqref{decomposition}, we obtain
\begin{align*}
u_0(0)&= A+\bA,\\
u'_0(0)&= \rmi \ell_c(A-\bA)+B+\bB,\\
u''_0(0)&=-\ell_c^2(A+\bA)+2\rmi \ell_c (B-\bB),\\
u'''_0(0)&=-\rmi \ell_c^3(A-\bA)-3 \ell_c^2 (B+\bB),
\end{align*}
from which we get a matrix passage from the quadruplet $(u_0(0),u'_0(0),u''_0(0),u'''_0(0))$ to $(A,\bA,B,\bB)$
\begin{equation*}
\mathscr{M}=\left(
\begin{matrix}
1 & 1 & 0 & 0 \\
\rmi \ell_c & -\rmi\ell_c & 1 & 1 \\
-\ell_c^2 & -\ell_c^2 & 2\rmi \ell_c & -2\rmi\ell_c \\
-\rmi\ell_c^3 & \rmi\ell_c^3 & -3\ell_c^2 & -3 \ell_c^2
\end{matrix}
\right).
\end{equation*}
One verifies that $\det \mathscr{M} = -16 \ell_c^4\neq0$ and computes
\begin{equation*}
\mathscr{M}^{-1}=\left(
\begin{matrix}
\frac{1}{2} & -\frac{3\rmi}{4\ell_c} & 0 & -\frac{\rmi}{4\ell_c^3} \\
\frac{1}{2} & \frac{3\rmi}{4\ell_c} & 0 & \frac{\rmi}{4\ell_c^3} \\
-\frac{\rmi\ell_c}{4} & -\frac{1}{4} & -\frac{\rmi}{4\ell_c} & -\frac{1}{4\ell_c^2}\\
\frac{\rmi\ell_c}{4} & -\frac{1}{4} & \frac{\rmi}{4\ell_c} & -\frac{1}{4\ell_c^2}
\end{matrix}
\right).
\end{equation*}
Let us introduce the map $\mathbf{q}:\C^4\rightarrow H_{-\eta}^4$ defined as $\mathbf{q}(z_1,z_2,z_3,z_4)=z_1 \zeta_0 + z_2 \overline{\zeta_0}+z_3 \zeta_1 + z_4 \overline{\zeta_1}$. We can then define the projection $\Q: H_{-\eta}^4 \rightarrow \E_0$ as
\bqq
\Q(u) = \mathbf{q}\left[\mathscr{M}^{-1}\left(u(0),u'(0),u''(0),u'''0) \right)^\mathbf{T} \right].
\label{projectionQnfe}
\eqq
Let us remark that the above definition gives
\begin{equation*}
\Q(u)=\left(\frac{u(0)}{2}-\frac{3\rmi u'(0)}{4\ell_c} -\frac{\rmi u'''(0)}{4\ell_c^3} \right) \zeta_0 + \left(-\frac{\rmi\ell_c u(0)}{4}  -\frac{u'(0)}{4}  -\frac{\rmi u''(0)}{4\ell_c}  -\frac{u'''(0)}{4\ell_c^2} \right)\zeta_1 + \text{ c.c.},
\end{equation*}
where c.c. stands for complex conjugate.

\paragraph{Center manifold theorem.} We can easily check that Hypothesis (H2$_\mu$) is satisfied as $\F$ is the composition of a pointwise operator and a convolution operator with exponential localization, where the pointwise operator is defined through the function $S$ which we suppose to be analytic in both arguments. As a consequence, we can apply the parameter-dependent center manifold with symmetries, to obtain the existence of neighborhoods $\U_u$, $\U_{\mu_c}$ of $(0,\mu_c)$ in $\E_0 \times (0,+\infty)$ and a map $\Psi \in \cC^k(\U_u\times \U_{\mu_c},\ker \Q)$ with $\Psi(0,0)=D_u\Psi(0,0)=0$, which commutes with $\mS_{1,2}$, and  such that for all $\mu\in \U_{\mu_c}$ the manifold
\begin{equation*}
\M_0(\mu-\mu_c):=\left\{ u_0+\Psi(u_0,\mu-\mu_c) ~|~ u_0 \in \U_u \right\} 
\end{equation*}
contains the set of all bounded solutions of \eqref{NFE}. From now on, we denote $\lambda := \mu-\mu_c$ and write,
\begin{equation*}
\Psi(u_0,\lambda) = \Psi(A,\bA,B,\bB,\lambda), \quad \text{ for } u_0 = A \zeta_0+\overline{A \zeta_0}+B \zeta_1+\overline{B \zeta_1}.
\end{equation*}
The fact that $\Psi$ should commute with $\mS_2$ implies that
\begin{equation*}
\mS_2 \Psi(A,\bA,B,\bB,\lambda) = \Psi (\mS_2\cdot (A,\bA,B,\bB),\lambda),
\end{equation*}
which yields
\begin{equation*}
- \Psi(A,\bA,B,\bB,\lambda) = \Psi (-A,-\bA,-B -\bB,\lambda).
\end{equation*}
Thus, there will not be any quadratic term in the Taylor expansion of $\Psi$. From now on, we write
\begin{equation*}
\Psi(A,\bA,B,\bB,\lambda)=\sum_{l_1,l_2,p_1,p_2,r>1}A^{l_1}\bA^{l_2}B^{p_1}\bB^{p_2}\lambda^r\Psi_{l_1,l_2,p_1,p_2,r},
\end{equation*}
the Taylor expansion of $\Psi$. Our next objective is to compute the lower order terms of this expansion.

\paragraph{Terms of order $\cO(\lambda A)$, $\cO(\lambda B)$, $\cO(\lambda \bA)$ and $\cO(\lambda \bB)$.} We start by computing the leading order terms in $\lambda$ in the above Taylor expansion of $\Psi$. For example, the function $\Psi_{1,0,0,0,1}$ is solution of the equation
\begin{equation*}
\T \Psi_{1,0,0,0,1} + \K \ast \zeta_0 = 0, \text{ with }  \Psi_{1,0,0,0,1} \in \ker \Q
\end{equation*}
We first note that $\K \ast \zeta_0 = \kappa_{0,1}\zeta_0$, such that one should look for solutions of the form
\begin{equation*}
\Psi_{1,0,0,0,1}(x) = \alpha_0 x^2 \zeta_0(x) + \psi_{1,0,0,0,1}(x), \text{ with }  \psi_{1,0,0,0,1} \in \E_0.
\end{equation*}
We then find that $\alpha_0$ should satisfy
\begin{equation*}
\alpha_0 \mu_c \kappa_{2,1} + \kappa_{0,1}=0, \text{ and } \alpha_0 = -\frac{\kappa_{0,1}^2}{ \kappa_{2,1}}.
\end{equation*}
Recall, that $\Psi_{1,0,0,0,1} \in \ker\Q$ and so $\Q\left( \Psi_{1,0,0,0,1}\right)=0$. We then write $ \psi_{1,0,0,0,1}=a_0 \zeta_0 + a_1\overline{\zeta_0}+b_0 \zeta_1 +b_1 \overline{\zeta_1}$ where the complex coefficients $(a_0,a_1,b_0,b_1)$ solve
\begin{equation*}
\alpha_0 \Q(x^2 \zeta_0(x)) +a_0 \zeta_0 (x)+ a_1\overline{\zeta_0}(x)+b_0 \zeta_1 (x)+b_1 \overline{\zeta_1}(x) =0,
\end{equation*}
as $\Q( \psi_{1,0,0,0,1})= \psi_{1,0,0,0,1}$. We find a set a four equations
\begin{equation*}
\left( \alpha_0 \frac{3}{2\ell_c^2}+a_0,-\alpha_0 \frac{3}{2\ell_c^2}+a_1,-\alpha_0 \frac{2 \rmi }{\ell_c}+b_0,-\alpha_0 \frac{\rmi }{\ell_c}+b_1\right)=(0,0,0,0),
\end{equation*}
where we used the fact that 
\begin{equation*}
\Q(x^2 \zeta_0(x))= \frac{3}{2\ell_c^2} \zeta_0(x)-  \frac{3}{2\ell_c^2} \overline{\zeta_0}(x)-\frac{2 \rmi }{\ell_c} \zeta_1(x) -\frac{ \rmi }{\ell_c}\overline{\zeta_1}(x).
\end{equation*}
As a consequence, we obtain
\bqq
\Psi_{1,0,0,0,1}(x) = \alpha_0 \left[\left( x^2+\frac{2\rmi }{\ell_c}x -\frac{3}{2\ell_c^2}\right)\rme^{\rmi \ell_c x} +\left(\frac{\rmi }{\ell_c}x+\frac{3}{2\ell_c^2}\right)\rme^{-\rmi \ell_c x}\right], \quad \forall x \in \R.
\label{Psi10001}
\eqq
Using the reflection symmetry $\mS_1$, we directly have that $\Psi_{0,1,0,0,1} = \mS_1 \Psi_{1,0,0,0,1} $. 

Let us now compute the function $\Psi_{0,0,1,0,1}$, associated to terms of the form $\lambda B$. It solves the equation
\begin{equation*}
\T \Psi_{0,0,1,0,1} +\K \ast \zeta_1=0, \text{ with } \Psi_{0,0,1,0,1}\in \ker \Q.
\end{equation*}
We first remark that $\K\ast \zeta_1 = \kappa_{0,1} \zeta_1$, so that we look for a solution of the form
\begin{equation*}
\Psi_{0,0,1,0,1}(x) = (\alpha_2 x^2 +\alpha_1 x) \zeta_1(x) +  \psi_{0,0,1,0,1}(x), \text{ with } \psi_{0,0,1,0,1} \in \E_0,
\end{equation*}
to get
\begin{align*}
-2\mu_c\alpha_2 \kappa_{2,1}+\kappa_{0,1}&=0,\\
\mu_c \alpha_1 \kappa_{2,1} - \mu_c \alpha_2 \kappa_{3,1}=0.
\end{align*}
From this, we deduce
\begin{equation*}
\alpha_2 = \frac{\kappa_{0,1}^2}{2 \kappa_{2,1}},\quad \text{ and } \quad \alpha_1 =  \frac{ \kappa_{3,1} \kappa_{0,1}^2}{2\kappa_{2,1}^2}.
\end{equation*}
Similarly, we recall that $\Psi_{0,0,1,0,1} \in \ker\Q$ and so $\Q\left( \Psi_{0,0,1,0,1}\right)=0$. We then write $ \psi_{0,0,1,0,1}=a_0 \zeta_0 + a_1\overline{\zeta_0}+b_0 \zeta_1 +b_1 \overline{\zeta_1}$ where the complex coefficients $(a_0,a_1,b_0,b_1)$ solve
\begin{equation*}
\Q\left(\left(\alpha_1x^2+\alpha_2 x^3\right) \zeta_0(x)\right) +a_0 \zeta_0(x) + a_1\overline{\zeta_0}(x)+b_0 \zeta_1(x) +b_1 \overline{\zeta_1}(x) =0.
\end{equation*}
We find that
\begin{align*}
a_0&= \frac{3\rmi}{2\ell_c^3}( \alpha_2+\rmi\ell_c \alpha_1),\\
a_1&=-\frac{3\rmi}{2\ell_c^3}( \alpha_2+\rmi\ell_c \alpha_1),\\
b_0&=\frac{1}{2\ell_c^2}\left( 4\rmi \ell_c\alpha_1 +3 \alpha_2\right),\\
b_1&=\frac{1}{2\ell_c^2}\left( 2\rmi \ell_c\alpha_1 +3 \alpha_2\right).
\end{align*}
As a consequence, we have for all $x\in\R$
\begin{align}
\Psi_{0,0,1,0,1}(x) &= \left(\alpha_2 x^3 +\alpha_1 x^2+\frac{4\rmi \ell_c\alpha_1 +3 \alpha_2}{2\ell_c^2}x +\frac{3\rmi \alpha_2-3\ell_c \alpha_1}{2\ell_c^3}\right) \rme^{\rmi \ell_c x} \nonumber\\
&\qquad +  \left(\frac{ 2\rmi \ell_c\alpha_1 +3 \alpha_2}{2\ell_c^2}x -\frac{3\rmi\alpha_2-3\ell_c \alpha_1}{2\ell_c^3}\right) \rme^{-\rmi \ell_c x}.
\label{Psi00101}
\end{align}
Using the reflection symmetry $\mS_1$, we directly have that $\Psi_{0,0,0,1,1} = - \mS_1 \Psi_{0,0,1,0,1} $. We next compute cubic coefficients in the Taylor expansion of $\Psi$.

\paragraph{Terms of order $\cO(A^2\bA)$ and $\cO(\bA^2A)$.} Using once more the symmetry $\mS_1$, we have that if $\Psi_{2,1,0,0,0}$ is known, then we have $\Psi_{1,2,0,0,0}=\mS_1\Psi_{2,1,0,0,0}$. We easily check that $\Psi_{2,1,0,0,0}$ solves
\begin{equation*}
0 = \T \Psi_{2,1,0,0,0} +\K * \zeta_0, \text{ with } \Psi_{2,1,0,0,0}\in \ker\Q,
\end{equation*}
which gives
\begin{equation*}
\Psi_{2,1,0,0,0}(x)=  -\alpha_0 x^2 \zeta_0(x) +\psi_{2,1,0,0,0}(x), \text{ with } \psi_{2,1,0,0,0} \in \E_0,
\end{equation*}
where $\alpha_0=-\kappa_{0,1}^2/ \kappa_{2,1}$. Computations similar to the ones for the term $\cO(\lambda A)$ lead to
\bqq
\Psi_{2,1,0,0,0}(x) = -\alpha_0 \left[\left( x^2+\frac{2\rmi }{\ell_c}x -\frac{3}{2\ell_c^2}\right)\rme^{\rmi \ell_c x} +\left(\frac{\rmi }{\ell_c}x+\frac{3}{2\ell_c^2}\right)\rme^{-\rmi \ell_c x}\right], \quad \forall x \in \R.
\label{Psi21000}
\eqq

\paragraph{The reduced vector field.} The reduced vector field will be of the form
\begin{subequations}
\begin{align}
\frac{\md A}{\md x}&= f_1(A,\bA,B,\bB,\lambda),\\
\frac{\md B}{\md x}&= f_2(A,\bA,B,\bB,\lambda),
\end{align}
\label{systAB}
\end{subequations}
together with the equations for the complex conjugates. Recall, that $f_1$ and $f_2$ are obtained by computing
\begin{equation*}
\frac{\md}{\md x} \Q\left(\Phi(u_0(\cdot +x))\right)|_{x=0}=(f_1, \overline{f_1}, f_2, \overline{f_2}).
\end{equation*}
Note that, slightly  abusing notation, we identify elements in $\E_0$ with their representation in the basis $\left\{\zeta_0,\overline{\zeta_0},\zeta_1,\overline{\zeta_1} \right\}$. We also remark that $\Phi (u_0) = u_0 + \Psi(u_0,\lambda)$, such that $\Q\left(\Phi(u_0(\cdot +x))\right)= \Q \left( u_0(\cdot+x)\right)  + \Q\left(\Psi(u_0(\cdot+x),\lambda \right)$ where 
\begin{equation*}
\frac{\md}{\md x} \Q \left( u_0(\cdot+x)\right)|_{x=0} =(\rmi \ell_c A+B, -\rmi \ell_c \bA+\bB,\rmi \ell_c B, -\rmi \ell_c \bB).
\end{equation*}
As a consequence, it remains to compute $\frac{\md}{\md x} \Q\left(\Psi(u_0(\cdot +x))\right)|_{x=0}$ only. On can check for example that, from the expression of $\Psi_{1,0,0,0,1}$ and $\Psi_{2,1,0,0,0}$ in \eqref{Psi10001} and \eqref{Psi00101} respectively that 
\begin{align*}
\frac{\md}{\md x} \Q\left(\Psi_{1,0,0,0,1}(\cdot +x)\right)|_{x=0}&=2\alpha_0 \left(\frac{\rmi}{\ell_c},-\frac{\rmi}{\ell_c}, 1,1\right),\text{ for terms of order } \cO(\lambda A),\\
\frac{\md}{\md x} \Q\left(\Psi_{0,0,1,0,1}(\cdot +x)\right)|_{x=0}&=2\frac{\ell_c\alpha_1-3\rmi \alpha_2}{\ell_c^2} \left(\frac{\rmi}{\ell_c},-\frac{\rmi}{\ell_c}, 1,1\right),\text{ for terms of order } \cO(\lambda B).
\end{align*}
We then find that the linear part of system \eqref{systAB} is given by
\begin{subequations}
\begin{align}
\frac{\md A}{\md x}&= \rmi \ell_c A+B+\frac{2\rmi \alpha_0}{\ell_c}\lambda \left( A +\bA\right) + \frac{2a_0}{\ell_c}\lambda \left(B-\bB\right) ,\\
\frac{\md B}{\md x}&= \rmi \ell_c B+2\alpha_0 \lambda \left( A +\bA\right) -2\rmi a_0 \lambda \left(B-\bB\right),
\end{align}
\label{ABlin}
\end{subequations}
where we set $a_0:=(3\alpha_2+\rmi \ell_c \alpha_1)/\ell_c^2 \in \R$. Following \cite[Lemma 2.4]{scheel-wu:14}, we know that there exists a smooth linear map $L(\lambda)$ such that for sufficiently small $\lambda$, the linear change of variables $(A,\bA,B,\bB)^\mathbf{T}=L(\lambda)(C,\bC,D,\bD)^\mathbf{T}$ transforms the linear system \eqref{ABlin} into the normal form
\begin{subequations}
\begin{align}
\frac{\md C}{\md x}&=\rmi \ell(\lambda) C+D,\\
\frac{\md D}{\md x}&=\alpha(\lambda) C+\rmi \ell(\lambda)D,
\end{align}
\label{ABlin2}
\end{subequations}
with complex conjugates, where we have set 
\begin{align*}
\ell(\lambda) &:=\frac{\rmi}{2}\sqrt{2\ell_c^2-4\ell_c a_0 \lambda +2\ell_c \sqrt{-4\ell_c a_0 \lambda +\ell_c^2 +8 \alpha_0 \lambda}},\\
\alpha(\lambda)&:=\ell_c a_0 \lambda -\frac{\ell_c^2}{2}+\frac{\ell_c}{2}\sqrt{-4\ell_c a_0 \lambda +\ell_c^2 +8 \alpha_0 \lambda}.
\end{align*}
Note that we have the expansions
\begin{equation*}
\ell(\lambda) =\ell_c +(\alpha_0-a_0)\lambda +\mathrm{o}(\lambda), \text{ and } \alpha(\lambda)=2 \alpha_0\lambda +\mathrm{o}(\lambda).
\end{equation*}
We are now going to apply a cubic transformation to our full system \eqref{systAB} for $\lambda =0$, that is to  
\begin{subequations}
\begin{align}
\frac{\md A}{\md x}&= \rmi \ell_c A+B + g_1(A,\bA,B,\bB),\\
\frac{\md B}{\md x}&= \rmi \ell_c B + g_2(A,\bA,B,\bB),
\end{align}
\label{systABcubic}
\end{subequations}
where we have set 
\begin{equation*}
g_{1,2}(A,\bA,B,\bB)= \sum_{n_1+n_2+n_3+n_4=3} A^{n_1}\bA^{n_2}B^{n_3}\bB^{n_4} \mathbf{g}^{1,2}_{n_1,n_2,n_3,n_4}.
\end{equation*}
First, we obtain the following expression for terms of order $ \cO(A|A|^2)$ which is given by computing 
\begin{equation*}
\frac{\md}{\md x} \Q\left(\Psi_{2,1,0,0,0}(\cdot +x)\right)|_{x=0}=-2\alpha_0 \left(\frac{\rmi}{\ell_c},-\frac{\rmi}{\ell_c}, 1,1\right).
\end{equation*}
Once again, following the strategy developed in \cite[Lemma 2.6]{scheel-wu:14}, one can find homogeneous polynomials of degree 3 denoted $(\mathcal{N}_1,\mathcal{N}_2)$ in the complex variables $(E, F, \bE, \bF)$, such that the change of variables
\begin{align*}
A &=E + \mathcal{N}_1(E, F, \bE, \bF),\\
B&=F + \mathcal{N}_2(E, F, \bE, \bF),
\end{align*}
is well-defined in a neighborhood of the origin and transforms the system \eqref{systABcubic} into the normal form
\begin{subequations}
\begin{align}
\frac{\md E}{\md x}&= \rmi \ell_c E+F +\mathcal{O}\left(\left( |E|+|F|\right)^5 \right),\\
\frac{\md F}{\md x}&= \rmi \ell_c F -2\alpha_0 E|E|^2+h_1 F|E|^2+h_2E\left( E\bF - \bE F \right)+\mathcal{O}\left(\left( |E|+|F|\right)^5 \right),
\end{align}
\label{systEF}
\end{subequations}
for tow complex constants $(h_1,h_2)\in \C$. As a consequence, applying our two change of variables and denoting with $(\widetilde{A},\overline{\widetilde{A}},\widetilde{B},\overline{\widetilde{B}})$ the new variables, we obtain the following system into normal form to leading order 
\begin{subequations}
\begin{align}
\frac{\md \widetilde{A}}{\md x}&=\rmi \ell(\lambda) \widetilde{A}+\widetilde{B},\\
\frac{\md \widetilde{B}}{\md x}&=\alpha(\lambda) \widetilde{A}+\rmi \ell(\lambda)\widetilde{B} -2\alpha_0 \widetilde{A}|\widetilde{A}|^2+h_1 \widetilde{B}|\widetilde{A}|^2+h_2\widetilde{A}\left( \widetilde{A}\overline{\widetilde{B}} - \overline{\widetilde{A}} \widetilde{B} \right).
\end{align}
\label{ABnormal}
\end{subequations}
The higher order terms in the normal form are of order
\begin{equation*}
|\lambda|\left( |\widetilde{A}|+|\widetilde{B}|\right)^3+\left(|\widetilde{A}|+|\widetilde{B}|\right)^5.
\end{equation*}
Next, we pass to corotating frame with respect to the normal form symmetry, 
\begin{equation*}
\widetilde{A}(x) =\rme^{\rmi \ell(\lambda) x} \widehat{A}(x) \text{ and } \widetilde{B}(x) = \rme^{\rmi \ell(\lambda) x} \widehat{B}(x),
\end{equation*}
to get at leading order
\begin{align*}
\frac{\md \widehat{A}}{\md x}&=\widehat{B},\\
\frac{\md \widehat{B}}{\md x}&=\alpha(\lambda) \widehat{A} -2\alpha_0 \widehat{A}|\widehat{A}|^2+h_1 \widehat{B}|\widehat{A}|^2+h_2\widehat{A}\left( \widehat{A}\overline{\widehat{B}} - \overline{\widehat{A}} \widehat{B} \right).
\end{align*}
We finally scale the equation, exhibiting leading order terms:
\begin{equation*}
\hat x = |\lambda|^{1/2}x, \quad \widehat{A} = |\lambda|^{1/2} \mathbf{A}, \quad  \widehat{B} = |\lambda| \mathbf{B},
\end{equation*}
which leads the new system
\begin{subequations}
\begin{align}
\frac{\md \mathbf{A}}{\md \hat x}&=\mathbf{B}+\mathcal{O}\left(|\lambda|^{1/2} \right),\\
\frac{\md \mathbf{B}}{\md \hat x}&=2\alpha_0 \mathbf{A}\left(\mathrm{sign}(\lambda)  -|\mathbf{A}|^2\right)+\mathcal{O}\left(|\lambda|^{1/2} \right).
\end{align}
\label{ABscale}
\end{subequations}

From now on, we assume that $\lambda>0$ and $\alpha_0>0$ which is equivalent to $\kappa_{2,1} <0$. In that case, we follow the perturbative analysis of \cite{iooss-per:93} (see also \cite{faye-etal:13}) and find a pair of reversible homoclinic orbits to the origin, solutions to \eqref{Nfe}, which can be approximated to leading order by
\begin{equation*}
u(x)=2\sqrt{\lambda} \mathrm{sech}\left(x\sqrt{2\alpha_0 \lambda} \right)\cos(x+\vartheta)+\cO(\lambda),\quad  x\in \R,
\end{equation*}
with $\vartheta\in\left\{0,\pi\right\}$ and $\lambda=\sqrt{\mu-\mu_c}$ for $\mu>\mu_c$.

\begin{rmk}
The example illustrates the somewhat novel (when compared to computations for local differential equations)  algebra involved with computing Taylor jets of the reduced vector field. We computed only the relevant cubic terms, that is, terms that give leading order expansions after scaling. Since the computation of those terms is somewhat simplified to a general computation of a reduced vector field, we include in the Appendix \ref{a:3rd} a computation of the vector field up to order 3.
\end{rmk}

\subsection{Slowly varying traveling waves in neural field equations}\label{s:tw}

Our second example is concerned with traveling waves rather than stationary solutions, in a system of $n$ coupled neural field equations,
\bqq
\partial_t\bu(t,x) = - D \bu(t,x)+\int_\R\rK(x-y) F(\bu(t,y),\mu)\md y, \quad (t,x) \in (0,\infty)\times\R,
\label{NFEn}
\eqq
for $\bu:\R\rightarrow \R^n$, $n\geq 1$, and $\mu\geq 0$, where $D=\text{diag}(d_j)$ is a diagonal matrix with positives entries $d_j>0$ for all $j=1\cdots n$. Throughout the sequel, we will assume that $\rK$ is a Gaussian matrix kernel in the sense that for all $1\leq i,j\leq n$, there exists $a_{i,j}>0$, such that $\rK_{i,j}(x) = \exp( - a_{i,j} x^2)$ for all $x\in\R$, and thus $\rK$ satisfies Hypothesis (H1) for all $\eta_0>0$. We also suppose that the nonlinear operator $\bu \mapsto \rK \ast F(\bu,\mu)$ verifies Hypothesis (H2$_\mu$) and that $\bu\mapsto F(\bu,\mu)$ is odd. Although this last assumption on the oddness of the nonlinearity is not required for the analysis and could be removed, it simplifies the subsequent computations of the reduced vector field on the center manifold. 

Spatially homogeneous states of \eqref{NFEn} are solutions of the kinetic equation on $\R^n$
\bqq
\frac{\md \bu }{\md t} = - D \bu + \rK_0 F(\bu,\mu),
\label{kinetic}
\eqq
where the matrix $\rK_0$ is defined through $ \rK_0:=\int_\R\rK(x)\md x$. In a neighborhood of $(\bu,\mu)=(0,0)$, we assume that the dynamics of \eqref{kinetic} can be reduced to a one-dimensional center manifold with a vector field
\begin{equation*}
\frac{\md z}{\md t} = g(z,\mu), \quad z \in\R.
\end{equation*}
We suppose that the resulting bifurcation is a supercritical pitchfork bifurcation.
\begin{hyp}[Supercritical pitchfork bifurcation]
The reduced vector field on the one-dimensional center manifold is odd in $z$ for all $\mu$ close to zero and
\begin{equation*}
g(z,\mu)=z\left( \alpha \mu - \beta z^2\right) +\mathcal{O}\left( |z|\left( |\mu|+z^2\right)^2\right), \text{ as } (z,\mu)\rightarrow (0,0)
\end{equation*}
with $\alpha>0$ and $\beta>0$.
\end{hyp}

Traveling wave solutions of \eqref{NFEn} are stationary solutions of the following system of equations
\bqq
\partial_t\bu = c \partial_\xi \bu - D \bu+\rK\ast  F (\bu,\mu),
\label{NFEnTW}
\eqq
where $\xi=x-ct$ for some constant $c\in \R$. Steady states of \eqref{NFEnTW} are thus solutions of the following nonlocal system
 \bqq
0 =\bu+\rG_c \ast  F (\bu,\mu),
\label{twn}
\eqq
where we set $\rG_c = \left(c\frac{\md}{\md \xi}I_n -D \right)^{-1}\rK$. It is important to note that $c \mapsto \rG_c$ is a smooth operator from $W^{1,\infty}(\R,\R^n)$ to itself because of the Gaussian nature of $\rK$. From now on, we will assume that there is a dependence between $c$ and $\mu$ by imposing that $c=\epsilon c_*$, $\mu=\epsilon^2$ for $\epsilon \geq 0$ and some $c_*\in\R$ independent of $\epsilon$. Such a scaling is motivated by an analogous study \cite{kirch-raugel:96} for systems of reaction-diffusion equations. It is also useful to note that in the limit $c\rightarrow 0$, we have $\rG_0=-D^{-1}\rK$.

The linearization of \eqref{twn} about the trivial state $\bu=0$ leads to the linear operator
\begin{equation*}
\T_{\epsilon} \bu:=\bu + \rG_{\epsilon c_*} \ast  \mathrm{D}_\bu F \left(0,\epsilon^2\right).
\end{equation*} 
We define the linear characteristic equation $d(\nu,\epsilon)$ as
\begin{equation*}
d(\nu,\epsilon ):=\det\left(\widehat{\T_\epsilon}(\nu) \right)=\det \left( I_n + \widehat{\rG_{\epsilon c_*}}(\nu)\mathrm{D}_\bu F \left(0,\epsilon^2\right) \right), \text{ for } (\nu,\epsilon)\in \C\times\R^+.
\end{equation*}
We make the following hypotheses on the characteristic equation.
\begin{hyp}[Homogeneous instability]
We assume that the characteristic equation $d(\nu,\epsilon)$ satisfies:
\begin{itemize}
\item $d(0,0)=\partial_\nu d(0,0)=0$ with $\partial_{\nu\nu}d(0,0)\neq 0$;
\item $d(i\ell,0)\neq0$ for all  $\ell \neq 0$.
\end{itemize}
\end{hyp}

\paragraph{Notation.} As $d(0,0)=0$, there exists $\be_0, \be_0^*\in\R^n$ such that
\begin{align*}
\widehat{\T_0}(0)\be_0 = \be_0+\widehat{\rG_{0}}(0) \mathrm{D}_\bu F(0,0)\be_0&=0,\\
\widehat{\T_0}(0)^{\mathbf{T}}\be_0^*=\be_0^*+\mathrm{D}_\bu F(0,0) ^{\mathbf{T}}\widehat{\rG_{0}}(0)^\mathbf{T}\be_0^*&=0,\\
\langle \be_0,\be_0^*\rangle&=1,
\end{align*} 
where $\langle\cdot,\cdot\rangle$ denotes the standard inner product on $\R^n$ given by
\begin{equation*}
\langle\bu,\bv\rangle=\sum_{k=1}^nu_kv_k, \text{ for any } \bu=(u_k)_{k=1}^n\in\R^n \text{ and } \bv=(v_k)_{k=1}^n\in\R^n.
\end{equation*}
Note that $\widehat{\rG_{0}}(0)=-D^{-1} \rK_0$, together with
\begin{align*}
\alpha &= \langle \rK_0 \mathrm{D}_{\bu,\mu}F(0,0)\be_0,\be_0^*\rangle>0,\\
\beta&= -\frac{1}{6}\langle \rK_0 \mathrm{D}_{\bu,\bu,\bu}F(0,0)\left[\be_0,\be_0,\be_0\right],\be_0^*\rangle>0,
\end{align*}
where $\alpha$ and $\beta$ are the coefficients appearing in the Taylor expansion of $g(z,\mu)$.

\paragraph{Symmetries.} As in the previous section, in addition to the translation equivariance, equation \eqref{twn} possesses two other symmetries, that we denote $\mS_1$ and $\mS_2$ respectively and act on functions as
\begin{equation*}
\mS_1 u(\xi) := u(-\xi), \quad \text{and }\quad \mS_2 u(\xi) :=-u(\xi), \quad \forall \xi \in \R.
\end{equation*}
The first symmetry is a consequence of the fact that each element of the matrix kernel $\rK$ is a symmetric function, whereas the second symmetry results from the odd symmetry of the nonlinear operator $F$ with respect to its first argument. Finally, let us remark that the conditions on the  dispersion relation ensures that the kernel $\E_0$ of $\T_0$ is given by
\begin{equation*}
\E_0 =\text{Span} \left\{ \be_0, \xi \be_0 \right\} \subset H^1_{-\eta}(\R,\R^n),
\end{equation*}
for all $0<\eta<\eta_0$ and any fixed $\eta_0>0$. As a consequence, any functions $\bu_0\in \E_0$, can be decomposed as
\bqq
\label{decn}
u_0 = A \be_0+B\be_1,
\eqq
for $(A,B)\in \R^2$ and $\be_1(\xi) := \xi \be_0$. We remark that the actions of $\mS_{1,2}$ on $\bu_0$ are given by 
\begin{align*}
\mS_1 u_0 &= A \be_0 - B\be_1, \\
\mS_2 u_0 &= -A \be_0- B\be_1.
\end{align*}
We identify the action of $\mS_{1,2}$ on the couple $(A,B)$ as
\begin{align*}
\mS_1 \cdot(A,B) &= (A, -B), \\
\mS_2 \cdot(A,B) &= (-A,-B).
\end{align*}

\paragraph{Projection $\Q$.} We now define the projection $\Q$ from $H^2_{-\eta}(\R,\R^n)\rightarrow \E_0$. Note that by Sobolev embedding we have $H^2(\R,\R^n)\subset \mathscr{C}^1(\R,\R^n)$, and thus we can take linear combinations of $\bu(0)$ and $\bu'(0)$. We define the projection $\Q:H^2_{-\eta}(\R,\R^n)\rightarrow \E_0$ through
\bqq
\Q(\bu):=\left( \bu(0),\be_0^*\right) \be_0 + \left( \bu'(0),\be_0^*\right)\be_1.
\label{projQn}
\eqq

\paragraph{Center manifold theorem.} We apply the parameter-dependent center manifold theorem with symmetries to system \eqref{twn}, to obtain the existence of neighborhoods $\U_\bu$, $\U_{0}$ of $(0,0)$ in $\E_0 \times (0,+\infty)$ and a map $\Psi \in \cC^k(\U_\bu\times \U_0,\ker \Q)$ with $\Psi(0,0)=D_\bu\Psi(0,0)=0$, which commutes with $\mS_{1,2}$, and  such that for all $\epsilon \in \U_0$ the manifold
\begin{equation*}
\M_0(\epsilon):=\left\{ \bu_0+\Psi(\bu_0,\epsilon) ~|~ \bu_0 \in \U_\bu \right\} 
\end{equation*}
contains the set of all bounded solutions of \eqref{twn}. From now on, we  write
\begin{equation*}
\Psi(\bu_0,\epsilon) = \Psi(A,B,\epsilon), \quad \text{ for } \bu_0 = A \be_0+B \be_1.
\end{equation*}
The fact that $\Psi$ should commute with $\mS_2$ implies that
\begin{equation*}
\mS_2 \Psi(A,B,\epsilon) = \Psi (\mS_2\cdot (A,B),\epsilon),
\end{equation*}
which yields
\begin{equation*}
- \Psi(A,B,\epsilon) = \Psi (-A,-B,\epsilon).
\end{equation*}
Thus, there will not be any quadratic term in the Taylor expansion of $\Psi$. From now on, we write
\begin{equation*}
\Psi(A,B,\epsilon)=\sum_{l_1,l_2,r>1}A^{l_1}B^{l_1}\epsilon^r\Psi_{l_1,l_2,r},
\end{equation*}
the Taylor expansion of $\Psi$. Our next task is to compute the lower order terms of this expansion.

\paragraph{Terms of order $\cO(\epsilon A)$ and $\cO(\epsilon B)$.} We first start by computing the linear leading order terms in $\epsilon$ in the above Taylor expansion of $\Psi$. The function $\Psi_{1,0,1}$ is the  solution to the equation
\begin{equation*}
\T_0 \Psi_{1,0,1}-c_*D^{-2}\frac{\md}{\md \xi} \left[\rK\ast \left(\mathrm{D}_\bu F(0,0)\be_0\right)\right]=0, \text{ with } \Psi_{1,0,1}\in \ker \Q.
\end{equation*}
A trivial computation shows that $\frac{\md}{\md \xi} \left[\rK\ast \left(\mathrm{D}_\bu F(0,0)\be_0\right)\right]=0$, such that $ \Psi_{1,0,1}\in \ker \Q \cap \ker \T_0$ and thus
\begin{equation*}
 \Psi_{1,0,1}=0.
\end{equation*}
On the other hand, we have that $\Psi_{0,1,1}$ is the solution to the equation
\begin{equation*}
\T_0 \Psi_{0,1,1}-c_*D^{-2}\frac{\md}{\md \xi} \left[\rK\ast \left(\mathrm{D}_\bu F(0,0)\be_1\right)\right]=0, \text{ with } \Psi_{0,1,1}\in \ker \Q.
\end{equation*}
First we note that, $\frac{\md}{\md \xi} \left[\rK\ast \left(\mathrm{D}_\bu F(0,0)\be_1\right)\right]=\rK_0 \mathrm{D}_\bu F(0,0)\be_0$ and we look for solutions of the form
\begin{equation*}
\Psi_{0,1,1}(\xi) = \gamma_0 \xi^2 \be_0 + \psi_{0,1,1},  \text{ with } \psi_{0,1,1}\in \E_0.
\end{equation*}
We then find that
\begin{equation*}
-\gamma_0 D^{-1}\int_\R y^2\rK(y)\mathrm{D}_\bu F(0,0)\be_0\md y - c_*D^{-2} \rK_0 \mathrm{D}_\bu F(0,0)\be_0=0,
\end{equation*}
such that
\begin{equation*}
\gamma_0=-\frac{c_*}{\kappa_2}, \quad \kappa_2:= \int_\R y^2\langle\rK(y)\mathrm{D}_\bu F(0,0)\be_0, \be_0^*\rangle\md y;
\end{equation*}
here, we used the fact that $\be_0=D^{-1}\rK_0 \mathrm{D}_\bu F(0,0) \be_0$ and $\langle\be_0,\be_0^*\rangle=1$. Note that $\kappa_2\neq 0$ as $\partial_{\nu\nu} d(0,0)\neq0$ from our hypothesis on the characteristic equation. Finally, as $\Q(\Psi_{0,1,1})=\psi_{0,1,1}$ and $\Psi_{0,1,1}\in\ker\Q$, we necessarily have $\psi_{0,1,1}=0$.

\paragraph{Terms of order $\cO(A^3)$.} The function $\Psi_{3,0,0}$ solves
\begin{equation*}
\T_0 \Psi_{3,0,0}-D^{-1}\rK\ast \left(\frac{1}{6}\mathrm{D}_{\bu,\bu,\bu} F(0,0)[\be_0,\be_0,\be_0]\right)=0, \text{ with } \Psi_{3,0,0}\in \ker \Q.
\end{equation*}
We find that
\begin{equation*}
\Psi_{3,0,0}(\xi) = \beta_0 \xi^2\be_0, 
\end{equation*}
where $\beta_0$ is given by
\begin{equation*}
\beta_0=-\frac{1}{6}\frac{\langle \rK_0 \mathrm{D}_{\bu,\bu,\bu} F(0,0)[\be_0,\be_0,\be_0],\be_0^*\rangle}{\kappa_2}=\frac{\beta}{\kappa_2}.
\end{equation*}

\paragraph{Terms of order $\cO(\epsilon^2 A)$.} The function $\Psi_{1,0,2}$ is solution of the equation
\begin{equation*}
\T_0 \Psi_{1,0,2}-D^{-1}\rK\ast \left(\mathrm{D}_{\bu,\mu} F(0,0)\be_0\right)=0, \text{ with } \Psi_{1,0,2}\in \ker \Q.
\end{equation*}
We find 
\begin{equation*}
\Psi_{1,0,2}(\xi) = \alpha_0 \xi^2\be_0, 
\end{equation*}
where $\alpha_0$ is given by
\begin{equation*}
\alpha_0=-\frac{\langle \rK_0 \mathrm{D}_{\bu,\mu}(0,0)\be_0,\be_0^*\rangle}{\kappa_2}=-\frac{\alpha}{\kappa_2}.
\end{equation*}

\paragraph{The reduced vector field.} The reduced vector field will be of the form
\begin{subequations}
\begin{align}
\frac{\md A}{\md \xi}&= f(A,B,\epsilon),\\
\frac{\md B}{\md \xi}&= g(A,B,\epsilon),
\end{align}
\label{systABn}
\end{subequations}
where $f$ and $g$ are obtained by computing 
\begin{equation*}
\frac{\md}{\md \xi} \Q\left(\Phi(\bu_0(\cdot +\xi))\right)|_{\xi=0}=(f,g).
\end{equation*}
Note that we slightly abused notation as we identify elements in $\E_0$ with their components on the basis $\left\{\be_0,\be_1 \right\}$. We also remark that $\Phi (\bu_0) = \bu_0 + \Psi(\bu_0,\epsilon)$, such that $\Q\left(\Phi(\bu_0(\cdot +\xi))\right)= \Q \left( \bu_0(\cdot+\xi)\right)  + \Q\left(\Psi(\bu_0(\cdot+\xi),\epsilon \right)$ where 
\begin{equation*}
\frac{\md}{\md \xi} \Q \left( \bu_0(\cdot+\xi)\right)|_{\xi=0} =(B,0).
\end{equation*}
Furthermore, we also have that
\begin{equation*}
\frac{\md}{\md \xi} \Q \left((\cdot+\xi)^2\be_0\right)|_{\xi=0} =(0,2).
\end{equation*}
Collecting all terms, we obtain the system
\begin{subequations}
\begin{align}
\frac{\md A}{\md \xi}&= B+\cO\left(\epsilon(|A|+|B|)+(|A|+|B|)^3\right),\\
\frac{\md B}{\md \xi}&= 2\epsilon \gamma_0 B +2 \alpha_0 \epsilon^2 A + 2 \beta_0 A^3+\cO\left( |B|(\epsilon^2 + |B|^2 +|A|^2)  \right).
\end{align}
\label{newABn}
\end{subequations}
We now rescale space with $\zeta= \epsilon \xi$, and the amplitudes $A=\epsilon \hat A$, $B = \epsilon^2 \hat B$ to obtain a new system
\begin{subequations}
\begin{align}
\frac{\md \hat A}{\md \zeta}&= \hat B+\cO\left(\epsilon \right),\\
\frac{\md \hat B}{\md \zeta}&= \frac{2}{\kappa_2} \left( -c_*  \hat B - \hat A  \left[\alpha - \beta \hat A^2\right] \right)+\cO\left( \epsilon \right).
\end{align}
\label{hatABn}
\end{subequations}

From now on we suppose that
\begin{equation*}
\kappa:=\frac{\kappa_2}{2}= \frac{1}{2}\int_\R y^2\langle\rK(y)\mathrm{D}_\bu F(0,0)\be_0, \be_0^*\rangle\md y>0,
\end{equation*}
and formally set $\epsilon =0$ in  \eqref{hatABn} to obtain the second order ordinary differential equation
\bqq
\kappa \frac{\md^2 \hat A}{\md \zeta^2} +c_* \frac{\md \hat A}{\md \zeta} + \hat A  \left[\alpha - \beta \hat A^2\right] =0.
\label{ODE}
\eqq
We know that such an equation admits monotone front solutions for any $|c_*|\geq 2\sqrt{\kappa \alpha}$ connecting the state $\hat A =0$ to the state $\hat A = \sqrt{\alpha/\beta}$ (see \cite{fisher:37,kpp:37}). Note that for $\epsilon>0$, there exists a unique saddle-point $\mathbf{a}(\epsilon):=(\sqrt{\alpha/\beta}+\cO(\epsilon),0)$. Then it follows from perturbative arguments \cite{eckmann-wayne:91,kirch-raugel:96} that system \eqref{hatABn} has front solutions connecting $(0,0)$ with $\mathbf{a}(\epsilon)$. Monotonicity in the tails can be established for speeds $|c_*|>2\sqrt{\kappa\alpha}+\cO(\epsilon)$. We denote by $u_*$ the front solution of equation \eqref{ODE} connecting $0$ to $\sqrt{\alpha/\beta}$, In our initial problem,  we thus have thus shown the existence of slowly varying front solutions of \eqref{NFEn} of the form
\begin{equation*}
\bu(t,x) = \epsilon u_*(\epsilon(x-\epsilon c_*t)) \be_0 +\cO(\epsilon^2),
\end{equation*}
for all $t,x\in\R$, with monotone tails for  $|c_*|\geq 2\sqrt{\kappa \alpha}+\cO(\epsilon)$ and $u_*$ solution of \eqref{ODE}.

\section{Discussion}

We established the existence of finite-dimensional center manifolds for nonlocal equations on the real line possessing a continuous translation symmetry. Rather than constructing a phase space, we use Lyapunov-Schmidt reduction on the set of trajectories to reduce to a finite-dimensional kernel, on which we construct a reduced flow with associated vector field through the shift induced by the action of translations on bounded solutions. There are clearly numerous generalizations possible, but also some apparent limitations to our approach, and we comment on a few of those here.

\paragraph{Nonautonomous systems.}

One can clearly allow for nonlinearities to depend on time explicitly $\F=\F(u(\cdot),\cdot)$. The reduction procedure remains literally unchanged. The explicit time dependence would however break the translation symmetry. One inherits an action of the shift mapping solutions $u(\cdot) \mapsto u(\cdot+\tau)$ combined with a shift in the time variable of the nonlinearity, $\F(u,\cdot)\mapsto \F(u,\cdot+\tau)$.  The reduced equations will inherit an equivalent action of the translation group, given by a non-autonomous vector field with similar time-dependence, for instance, periodic, quasi-periodic, heteroclinic, etc. In this light, our approach can be viewed as an analysis similar to the constructions of trajectory attractors in non-autonomous or ill-posed evolution equations; see for instance \cite{sell,vis}.

\paragraph{Infinite-dimensional systems.}

We studied nonlocal equations where $u(x)\in\R^n$. It would be interesting and quite useful to generalize to equations where $u(x)\in \X$, a Hilbert space or even a Banach space. The main obstacle at this point is the fact that the results in \cite{faye-scheel:13} are limited to finite-dimensional ranges. It is conceivable that those results could be generalized with suitable compactness assumptions on lower-order terms.

\paragraph{Semilinear equations only.}

Another limitation of our results, again owed to the limitations in \cite{faye-scheel:13}, is the fact that we require our equations to be semilinear in the sense that the principal part in the sense of regularity is invertible, chosen as the identity, here, and other terms are of lower order, somewhat regular convolution kernels. Results in \cite{afss} motivate that significantly different phenomena can be expected when such hypotheses are violated. In particular, one may find non-smooth solutions, precluding the possibility of a differentiable action of the translation group. Examples are in particular kernels containing Dirac-masses, such as kernels mimicking lattice differential equations $\K(\cdot)=\sum a_j\delta(\cdot-\xi_j)$.

\paragraph{Exponentially  localized kernels only.}

We rely on exponential localization of the kernel when invoking \cite{faye-scheel:13}, and also when formulating our contraction-mapping theorem in spaces of exponentially growing functions. It is not clear how to weaken those assumptions significantly. It is conceivable to formulate assumptions on multiplicities of roots of the characteristic equation $d(\nu)$ on the imaginary axis given sufficient algebraic localization, such that moments of $\K$ are well defined and sufficiently high moments do not vanish, but Fredholm properties of the linear part as well as the nonlinear arguments will likely require different choices of function spaces and possibly additional assumptions on the nonlinearity. In this direction, establishing Fredholm properties and henceforth existence of center manifolds in exponentially weighted spaces of continuous functions, as used in \cite{vander-vangils:87}, would yield sharper results since  nonlinearities $f\in \cC^k$ yield superposition operators $\mathcal{F}$ of class $\cC^k$.

\paragraph{Center-manifolds versus asymptotic methods.}

Without reviewing the general merits of center manifolds, we would like to point out that center manifold methods have inherent advantages compared to more direct matched asymptotics or scaling arguments. One could for instance try to find solutions on the center manifold directly, using formally derived leading-order approximations, and control errors in a subsequent step, locally near a specific solution. While such approaches may be more robust, for instance in the case of algebraic localization of the kernel, they give weaker results; see for instance \cite{fs}. In particular, uniqueness statements are restricted to neighborhoods of particular ansatz solutions. Statements such as a Poincar\'e-Bendixson theorem, immediate for two-dimensional kernels, seem elusive without the construction of an actual flow.

\paragraph{Nonlocal equations versus local ODEs.}

Inspecting the calculations of the reduced vector field, one realizes that only finitely many generalized moments \eqref{eq:kappa} enter the computation of Taylor jets at any fixed order. One can therefore formally replace the nonlocal convolution kernel with a differential equations that reflects the Taylor jet of the Fourier transform of the convolution kernel at roots of the characteristic equation $d(\nu)=0$, and find the exact same reduced differential equation. Since higher orders of the Taylor jet invoke higher generalized moment, there may however not be one differential equation that yields the vector field associated with a fixed given nonlocal equation. Computationally, this observation does however provide an alternative strategy towards computing the reduced vector field, relying on an approximating differential equation and the more tradition computation of center manifolds in the associated phase space. It would be interesting to formulate a priori conditions for sufficiently high orders of approximation.

\appendix
\section{Superposition operators on exponentially weighted spaces}\label{a:non}

We show how to adapt  \cite[Lemma 3 \& 5]{vander-vangils:87} in order to show that Hypothesis (H2) holds in the case of a nonlinearity 
\begin{equation*}
\F(u)[x]=g(u(x)), \quad \forall x \in \R, \quad \forall u \in W^{1,\infty}(\R,\R^n),
\end{equation*}
with $g\in \cC^{k+1}$. For convenience of the presentation, we simply denote $H^1_{-\eta}$ instead of $H^1_{-\eta}(\R,\R^n)$ throughout this section.

Let us first suppose first that $g\in \mathscr{C}_b^1(\R^n)$, the set of $\cC^1$ bounded functions, and write $g_0=\sup |g|$, $g_1=\sup |g'|$. For $\eta \geq \zeta>0$ let us show that $\F:H^1_{-\zeta}\rightarrow H^1_{-\eta}$ is continuous. For this, let $u,v \in H^1_{-\zeta}$, and decompose
\[
\| \F(u)-\F(v)\|_{H^1_{-\zeta}}^2=\| \omega_{-\eta}\left(\F(u)-\F(v)\right)\|_{L^2}^2+\| \omega_{-\eta} \partial_x(\F(u)-\F(v))\|_{L^2}^2.
\]
We estimate the first integral as
\[
\| \omega_{-\eta}\left(\F(u)-\F(v)\right)\|_{L^2}^2\leq 2g_0^2 \int_{|x|\geq \rho} \omega_{-\eta}^2(x)\md x + \underset{|x|\leq \rho}{\sup} |g(u(x))-g(v(x))|^2  \| \omega_{-\eta}\|_{L^2}^2.
\]
Now, fix $\epsilon>0$ such that one can find some $\rho>0$ so that
\[
2g_0^2 \int_{|x|\geq \rho} \omega_{-\eta}^2(x)\md x \leq \frac{\epsilon^2}{4},
\]
Furthermore, since $\omega_{-\eta} u \in H^1$, by Morrey's inequality, we have that $\omega_{-\eta} u \in \mathscr{C}^{0,1/2}$, which in turn implies that the set $\Omega:=\left\{u(x) ~|~|x| \leq \rho \right\}$ is compact. As a consequence, since $g$ is continuous, there exists $\delta>0$ such that
\[
|g(y+z)-g(y)| < \frac{\epsilon}{2\| \omega_{-\eta}\|_{L^2}}, \text{ if } y \in \Omega, \text{ and } |z|< \delta_1.
\]
Let $\delta>0$, to be chosen later, and $u,v \in H^1_{-\zeta}$ with $\|u-v\|_{ H^1_{-\zeta}}\leq \delta$. Then, again by Morrey's inequality, we have that for some constant $C>0$,
\[
\| \omega_{-\zeta} (u-v)\|_{\mathscr{C}^{0,1/2}} \leq C\|u-v\|_{ H^1_{-\zeta}}\leq C\delta.
\]
Thus, for all $|x|\leq \rho$, 
\[
|u(x)-v(x)| \leq C \delta \underset{|x|\leq \rho}{\sup}\omega_{-\zeta}(x)^{-1}.
\]
As a consequence, we choose $\delta:= \left( C  \underset{|x|\leq \rho}{\sup}\omega_{-\zeta}(x)^{-1}\right)^{-1}$, and obtain
\[
\underset{|x|\leq \rho}{\sup} |g(u(x))-g(v(x))|^2 \| \omega_{-\eta}\|_{L^2}^2 \leq \frac{\epsilon^2}{4},
\]
and then
\[
\| \omega_{-\eta}\left(\F(u)-\F(v)\right)\|_{L^2}^2\leq \frac{\epsilon^2}{2}.
\]
Finally, we remark that
\[
\| \omega_{-\eta} \partial_x(\F(u)-\F(v))\|_{L^2}^2 \leq \int_\R \omega_{-\eta}^2(x) |u'(x)-v'(x)|^2|g'(u(x))|^2\md x + \int_\R \omega_{-\eta}^2(x) |v'(x)|^2|g'(u(x))-g'(v(x))|^2\md x.
\]
Note that the first integral is controlled by
\[
\int_\R \omega_{-\eta}^2(x) |u'(x)-v'(x)|^2|g'(u(x))|^2\md x \leq g_1^2 \|u-v\|_{ H^1_{-\zeta}}.
\]
The second integral can be evaluated similarly, using the fact that $g'$ is continuous and $\omega_{-\zeta} v' \in L^2$, such that we also get
\[
\| \omega_{-\eta} \partial_x(\F(u)-\F(v))\|_{L^2}^2 \leq \frac{\epsilon^2}{2}.
\]

We now turn to differentiability of the superposition operator. 
Let $g\in \mathscr{C}_b^{k+1}(\R^n)$, and let $\eta > k\zeta>0$ for $k\geq1$, then $\F:H^1_{-\zeta}\rightarrow H^1_{-\eta}$ is $\mathscr{C}^k$. First,  define for each $1\leq p \leq k$ a mapping $\F^{(p)}$ by $\F^{(p)}(u)(x):=D^pg(u(x))$ for any $x\in\R$ and $u\in W^{1,\infty}$. We consider $\F^{(p)}$ as a $p$-linear operator given through
\[
\F^{(p)}(u)\cdot(u_1,\ldots,u_p)(x):=D^pg(u(x))\cdot(u_1(x),\ldots,u_p(x)), \quad \forall x\in\R, \forall u_1,\ldots,u_p\in W^{1,\infty}.
\]
It is easy to check that $\F^{(p)}(u) \in \mathscr{L}^{(p)}(H^1_{-\zeta}\times \cdots \times H^1_{-\zeta},H^1_{-\eta})$. Indeed, from its definition, we have that
\[
\| \F^{(p)}(u)\cdot(u_1,\ldots,u_p) \|_{H^1_{-\eta}} \leq \| \F^{(p)}(u) \|_{H^1_{-\eta+p \zeta}} \| u_1 \|_{H^1_{-\zeta}} \cdots \| u_p \| _{H^1_{-\zeta}}, 
\]
and, as consequence, we also have that $\F^{(p)}$ is continuous from $H^1_{-\zeta}$ into $\mathscr{L}^{(p)}(H^1_{-\zeta}\times \cdots \times H^1_{-\zeta},H^1_{-\eta})$. Furthermore, for any $u,v \in H^1_{-\zeta}$, we have 
\begin{align*}
\| \F(u+v)-\F(u)-\F^{(1)}(u)\cdot v \| _{H^1_{-\eta}} &= \left\| \int_0^1 \left(\F^{(1)}(u+sv) -\F^{(1)}(u) \right)\cdot v \md s \right\|_{H^1_{-\eta}}\\
&\leq \underset{\ s\in [0,1]}{\sup}\| \F^{(1)}(u+sv) -\F^{(1)}(u)  \|_{H^1_{-\eta+p\zeta}} \| v \|_{H^1_{-\zeta}}.
\end{align*}
Since, in particular, $\F^{(1)}$ is continuous, we have for each $\epsilon>0$ the existence of $\delta>0$ such that
\[
\underset{\ s\in [0,1]}{\sup}\| \F^{(1)}(u+sv) -\F^{(1)}(u)  \|_{H^1_{-\eta+p\zeta}} \leq \epsilon, \text{ if } \| v \| _{H^1_{-\zeta}} \leq \delta.
\]
Thus, $\F$ is differentiable at $u \in H^1_{-\zeta}$, and one can prove in the same fashion that $\F$ is $\mathscr{C}^p$ for each $p=1,\ldots,k$.

% \item 
Finally, let $g\in \mathscr{C}^{k+1}(\R^n)$ and let $\chi_\epsilon$ be the cut-off operator introduced in the previous section, we define
\[
g^\epsilon(u):=g\left(\bar\chi(u/\epsilon) u\right),
\]
and, since $\bar\chi$ is chosen as a smooth cut-off, we have that $g^\epsilon \in \mathscr{C}^{k+1}_b(\R^n)$ and thus $\F^\epsilon:H^1_{-\zeta}\rightarrow H^1_{-\eta}$ is $\mathscr{C}^k$.

\section{Contractions on scales of embedded Banach spaces}\label{secapp:thm}

We give details on the proof of Lemma \ref{l:cmfdsm}. 

We first recall a result from \cite{vander-vangils:87} on contractions on embedded Banach spaces. Let $\X,\Y,\cZ$ and $\Lambda$ be Banach spaces with norms denoted respectively by $\| \cdot \|_\X$, $\| \cdot \|_\Y$, $\| \cdot \|_\cZ$ and $\| \cdot \|_\Lambda$, with continuous embedding:
\begin{equation*}
\X \overset{\J}{\hookrightarrow} \Y \overset{\G}{\hookrightarrow} \cZ.
\end{equation*}
Consider the fixed point equation
\bqq
\label{eqFP}
y=\mf(y,\lambda),
\eqq
where $\mf:\Y\times \Lambda \longrightarrow \Y$ satisfies the following conditions:
\begin{itemize}
\item[(C1)] $\G \mf:\Y\times \Lambda \longrightarrow \cZ$ has continuous partial derivative $D_y(\G \mf): \Y\times \Lambda \longrightarrow \cL(\Y,\cZ)$ with
\begin{equation*}
D_y(\G \mf)(y,\lambda)=\G \mf^{(1)}(y,\lambda)=\mf^{(1)}_1(y,\lambda) \G, \quad \forall (y,\lambda)\in \Y\times \Lambda,
\end{equation*}
for some $\mf^{(1)}:\Y\times \Lambda \longrightarrow \cL(\Y)$ and $\mf^{(1)}_1:\Y\times \Lambda \longrightarrow \cL(\cZ)$.
\item[(C2)] $\mf_0:\X\times \Lambda \longrightarrow \Y$, $(y_0,\lambda)\longmapsto \mf_0(y_0,\lambda)=\mf(\J y_0,\lambda)$ has continuous partial derivative $D_\lambda \mf_0:\X\times \Lambda \longrightarrow \cL(\Lambda,\Y)$.
\item[(C3)] There exists $\kappa\in [0,1)$ such that 
\begin{equation*}
\| \mf(y,\lambda)-\mf(\tilde y, \lambda) \|_\Y\leq \|y-\tilde y \|_\Y, \quad \forall y,\tilde y \in \Y, \quad \forall \lambda \in \Lambda,
\end{equation*}
and
\begin{equation*}
\| \mf^{(1)}(y,\lambda)\|_\Y \leq \kappa, \quad \| \mf^{(1)}_1(y,\lambda) \|_\cZ \leq \kappa, \quad \forall (y,\lambda)\in \Y\times \Lambda.
\end{equation*}
\item[(C4)] Let $y=\tilde y(\lambda)\in \Y$ be the unique solution of \eqref{eqFP} for $\lambda \in \Lambda$. Suppose that $\tilde y(\lambda)= \J \tilde y_0(\lambda)$ for some continuous $\tilde y_0:\Lambda \longrightarrow \X$.
\end{itemize}

These conditions allow to consider the following equation in $\cL(\Lambda,\Y)$:
\bqq
\label{eqFPLin}
\Theta=\mf^{(1)}(\tilde y(\lambda),\lambda)\Theta+D_\lambda \mf_0(\tilde y_0(\lambda),\lambda),
\eqq
which has a unique solution $\tilde \Theta(\lambda)\in \cL(\Lambda,\Y)$ for any $\lambda\in \Lambda$ from condition $(C3)$. The following Theorem is proved in \cite{vander-vangils:87}.
\begin{thm}\label{thm:appendix}
Assume $(C1)-(C4)$. Then the solution map $\tilde y: \Lambda\rightarrow \Y$ of \eqref{eqFP} is Lipschitz continuous, and $\G \tilde y:\Lambda \rightarrow \cZ$ is of class $\cC^1$, with
\bqq
\label{eqResApp}
D_\lambda \G \tilde y(\lambda)=\G \tilde \Theta(\lambda), \quad \forall \lambda \in \Lambda.
\eqq
\end{thm}

We now turn to the proof of Lemma \ref{l:cmfdsm}, considering first the case $p=1$ and then higher-order differentiability. Once again, for convenience of the presentation, we simply denote $H^1_{-\eta}$ instead of $H^1_{-\eta}(\R,\R^n)$.

\paragraph{Continuous differentiability of the fixed point, $p=1$.} We fix $ \eta \in (\tilde\eta,\bar\eta]$ and apply Theorem \ref{thm:appendix} with $\X=\Y=H^{1}_{-\tilde\eta }$, $\cZ=H^{1}_{-\eta }$, $\Lambda=\E_0$ and $\mf(y,\lambda):=\cS^\epsilon(y;\lambda)$, where $\cS^\epsilon:H^{1}_{-\tilde\eta }\times \E_0\longrightarrow H^{1}_{-\tilde\eta }$ is defined in \eqref{eqNLeps}. Indeed, one can check that all assumptions (C1)-(C4) are met in that case, and we obtain that $\Phi:\E_0\rightarrow H^1_{-\eta}$ is of class $\cC^1$ with derivative $\Phi^{(1)}(u_0):=D\Phi(u_0)\in \mathscr{L}(\E_0, H^{1}_{-\eta })$ being the unique solution of the equation
\bqq
\Theta = D_y \mf(\Phi(u_0),u_0) \Theta + D_\lambda \mf(\Phi(u_0),u_0):=F_1(\Theta,u_0).
\label{eq:fixed_point}
\eqq
Note that the mapping $F_1: \mathscr{L}(\E_0, H^{1}_{-\eta }) \times \E_0\rightarrow  \mathscr{L}(\E_0, H^{1}_{-\eta })$ is a uniform contraction for each $\eta\in[\tilde\eta,\bar\eta]$ by the assumptions on $\F^\epsilon$ (recall that $D\F^\epsilon(u)$ is assumed to be Lipschitz in $u$), and hence we have that its fixed point $\Phi^{(1)}(u_0)$ belongs in fact to $\mathscr{L}(\E_0, H^{1}_{-\tilde\eta })$ by continuous embedding and thus $\Phi^{(1)}:\E_0\rightarrow \mathscr{L}(\E_0, H^{1}_{-\eta })$ is continuous if $\eta \in (\tilde\eta,\bar\eta]$.
\paragraph{Higher smoothness, $p\geq 2$.} We now use induction on $p$. Let $1\leq p < k$, and suppose that for all $q$ with $1\leq q \leq p$ and for all $\eta \in (q\tilde\eta,\bar\eta]$ the mapping $\Phi:\E_0\rightarrow H^1_{-\eta}$ is of class $\cC^p$, with $\Phi^{(q)}(u_0):=D^q\Phi(u_0)\in \mathscr{L}^{(q)}(\E_0, H^{1}_{-q\tilde\eta })$ for each $u_0\in \E_0$ and $\Phi^{(q)}: \E_0\rightarrow  \mathscr{L}^{(q)}(\E_0, H^{1}_{-\eta })$ continuous if $\eta \in (q\tilde\eta,\bar\eta]$. Suppose also that $\Phi^{(p)}(u_0)$ is the unique solution of an equation that is of the form
\bqq
\Theta^{(p)} = D_y \mf(\Phi(u_0),u_0) \Theta^{(p)}  + G_p(u_0):=F_p(\Theta^{(p)},u_0),
\label{eq:fixed_point2}
\eqq
with $G_1(u_0)=D_\lambda \mf(\Phi(u_0),u_0)$ and, for $p\geq 2$, $G_p(u_0)$ is given as a finite sum of terms of the form
\[
D^{(k)}_yD^{(q-k)}_\lambda\mf(\Phi(u_0),u_0)\cdot\left(D^{(r_1)}\Phi(u_0),\ldots,D^{(r_k)}\Phi(u_0) \right),
\]
with $2\leq q\leq p$ and $k \leq q$ with $1\leq r_i <p$ for all $i=1,\ldots,k$ that verify $r_1+\cdots+r_k=k$. We remark that we have $G_p(u_0)\in \mathscr{L}^{(p)}(\E_0, H^{1}_{-p\tilde\eta })$. As a consequence, the mapping $F_p:\mathscr{L}^{(p)}(\E_0, H^{1}_{-\eta })\times \E_0 \rightarrow \mathscr{L}^{(p)}(\E_0, H^{1}_{-\eta })$ is well defined and is a uniform contraction for all $\eta \in [p \tilde\eta,\bar\eta]$. However, the term $D_y \mf(\Phi(u_0),u_0)$ is not continuously differentiable and one needs to apply Theorem \ref{thm:appendix} using three different  Banach spaces. Therefore, fix some $\eta \in ((p+1)\tilde\eta,\bar\eta]$ and choose $\sigma \in (\tilde \eta, \eta /(p+1))$ and $\zeta \in ((p+1)\sigma,\eta)$. We now show that the hypotheses of the theorem are satisfied with $\X=\mathscr{L}^{(p)}(\E_0, H^{1}_{-p\sigma })$, $\Y=\mathscr{L}^{(p)}(\E_0, H^{1}_{-\zeta })$ and $\cZ=\mathscr{L}^{(p)}(\E_0, H^{1}_{-\eta })$, $\Lambda = \E_0$ and $\mf = F_p$. Condition (C3) is met since $C(\eta)\delta_1(\epsilon)<1$ for all $\eta \in [\tilde \eta,\bar \eta]$, while (C4) follows from the induction hypothesis and the fact that $\sigma > \tilde \eta$. One can then check that $u_0\mapsto D_y \mf(\Phi(u_0),u_0)$ is continuous from $\E_0$ into $\mathscr{L}(H^{1}_{-\zeta}, H^{1}_{-\eta })$ as $\eta>\zeta$ and that $\Phi:\E_0\rightarrow H^{1}_{-\zeta }$ is continuous. In fact, we further have that $u_0\mapsto D_y \mf(\Phi(u_0),u_0)$ is $\cC^1$ from $\E_0$ into $\mathscr{L}(H^{1}_{-p \sigma }, H^{1}_{-\zeta })$ which follows from the fact that $\zeta>(p+1)\sigma$ and that $\Phi:\E_0\rightarrow H^{1}_{-\sigma }$ is of class $\cC^1$. Provided that $G_p:\E_0\rightarrow H^{1}_{-\zeta }$ is of class $\cC^1$  we then conclude from Theorem \ref{thm:appendix} that $\Phi^{(p)}:\E_0\rightarrow \mathscr{L}^{(p)}(\E_0, H^{1}_{-\eta })$ is of class $\cC^1$ and hence $\Phi:\E_0\rightarrow H^{1}_{-\eta }$ if of class $\cC^{p+1}$ if $\eta \in ((p+1)\tilde\eta,\bar\eta]$. The proof of the fact that $G_p:\E_0\rightarrow H^{1}_{-\zeta }$ is of class $\cC^1$ follows along similar lines as \cite[Lemma 7]{vander-vangils:87} and is omitted here.

\section{Computations of cubic order terms of the Taylor expansion}\label{a:3rd}

We compute expansions to order 3 of the reduction function $\Psi$ and the reduced vector field for the example from Section \ref{s:stat}. 

\subsection{Expansion of the reduction function}
We calculate the general cubic terms of $\Psi$. 
\paragraph{Terms of order $\cO(A^3)$ and $\cO(\bA^3)$.} We first note that if $\Psi_{3,0,0,0,0}$ is know, then we have $\Psi_{0,3,0,0,0}=\mS_1\Psi_{3,0,0,0,0}$. collecting terms of order $\cO(A^3)$, we obtain the equation
\begin{equation*}
0 = \T \Psi_{3,0,0,0,0} -\frac{1}{3}\K * \zeta_0^3, \text{ with } \Psi_{3,0,0,0,0}\in \ker\Q.
\end{equation*}
By noting that $\K * \zeta_0^3 = \kappa_{0,3} \zeta_0^3$, we obtain that 
\begin{equation*}
\Psi_{3,0,0,0,0}(x)= \frac{1}{3}\frac{ \kappa_{0,3}}{-1+ \kappa_{0,3}/\kappa_{0,1}}e^{3\rmi\ell_c x}+\psi_{3,0,0,0,0}(x) \text{ with } \psi_{3,0,0,0,0} \in \E_0.
\end{equation*}
Then, one computes that $\Q(e^{3\rmi\ell_c x})=-4\zeta_0(x)+5\overline{\zeta_0}(x)+8\rmi \ell_c \zeta_1(x)+4\rmi \ell_c \overline{\zeta_1}(x)$, such that we have
\bqq
\Psi_{3,0,0,0,0}(x)= \frac{1}{3}\frac{ \kappa_{0,3}}{-1+ \kappa_{0,3}/\kappa_{0,1}}\left( e^{3\rmi\ell_c x} +\left(4-8\rmi \ell_cx\right)\rme^{\rmi\ell_c x}-(5+4\rmi\ell_cx)\rme^{-\rmi\ell_c x} \right), \quad \forall x \in \R.
\label{Psi3000}
\eqq

\paragraph{Terms of order $\cO(B^3)$ and $\cO(\bB^3)$.} By symmetry we have that $\Psi_{0,0,0,3,0}=-\mS_1\Psi_{0,0,3,0,0}$ where $\Psi_{0,0,3,0,0}$ solves
\begin{equation*}
0 = \T \Psi_{0,0,3,0,0} -\frac{1}{3}\K * \zeta_1^3, \text{ with } \Psi_{0,0,3,0,0} \in \ker\Q.
\end{equation*}
We first note that
\begin{equation*}
\K * \zeta_1^3 (x) = \left[ \kappa_{0,3} x^3 -2\kappa_{1,3}x^2 -2\kappa_{2,3}x - \kappa_{3,3}\right]\rme^{3\rmi \ell_c x}.
\end{equation*}
As a consequence, we look for solutions of the form
\begin{equation*}
 \Psi_{0,0,3,0,0} (x) = \left(\beta_0+\beta_1 x + \beta_2 x^2 + \beta_3  x^3 \right) e^{3\rmi \ell_c x} +\psi_{0,0,3,0,0}(x), \text{ with } \psi_{0,0,3,0,0} \in \E_0.
\end{equation*}
Collecting terms of same order, we find a recursive system of equations for $(\beta_j)_{j=0,\cdots,3}$
\begin{align*}
\beta_3\left( -1+\mu_c \kappa_{0,3} \right) -\frac{\kappa_{0,3}}{3}&=0,\\
\beta_2\left( -1+\mu_c \kappa_{0,3} \right) +\left(\frac{2}{3}-2\mu_c\beta_3\right) \kappa_{1,3} &=0,\\
\beta_1\left( -1+\mu_c \kappa_{0,3} \right) + \left(\frac{2}{3}-2\mu_c\beta_3\right)\kappa_{2,3} -2\mu_c \beta_2 \kappa_{1,3} &=0,\\
\beta_0\left( -1+\mu_c \kappa_{0,3} \right) +\left(\frac{1}{3}-\mu_c\beta_3\right)\kappa_{3,3}+\mu_c \beta_2 \kappa_{2,3}-\mu_c\beta_1 \kappa_{1,3} &=0.
\end{align*}
We find that
\begin{equation*}
\beta_3=\frac{1}{3}\frac{\kappa_{0,3}}{-1+ \kappa_{0,3}/\kappa_{0,1}},\quad \beta_2=\frac{2}{3}\frac{\kappa_{1,3}}{(-1+ \kappa_{0,3}/\kappa_{0,1})^2}, \quad \beta_1=\frac{2}{3}\frac{\kappa_{2,3}}{(-1+ \kappa_{0,3}/\kappa_{0,1})^2}+\frac{4}{3}\frac{\kappa_{1,3}^2/\kappa_{0,1}}{(-1+ \kappa_{0,3}/\kappa_{0,1})^3},
\end{equation*}
and
\begin{equation*}
\beta_0=\frac{1}{3}\frac{\kappa_{3,3}}{(-1+ \kappa_{0,3}/\kappa_{0,1})^2}+\frac{4}{3}\frac{\kappa_{1,3}^3/\kappa_{0,1}^2}{(-1+ \kappa_{0,3}/\kappa_{0,1})^4}.
\end{equation*}
straightforward computations show that
\begin{align*}
\psi_{0,0,3,0,0}&= - \frac{-8\ell_c^3 \beta_0 +9\ell_c \beta_2 +\rmi(12\ell_c^2\beta_1-3\beta_3)}{2\ell_c^3}\zeta_0+\frac{-10\ell_c^3 \beta_0 +9\ell_c \beta_2 +\rmi(12\ell_c^2\beta_1-3\beta_3)}{2\ell_c^3}\overline{\zeta_0}\\
&\qquad  - \frac{16\ell_c^2\beta_1-3\beta_3 +\rmi(16 \ell_c^3 \beta_0 - 10 \ell_c \beta_2)}{2\ell_c^2}\zeta_1-\frac{10\ell_c^2\beta_1-3\beta_3 +\rmi(8 \ell_c^3 \beta_0 - 8 \ell_c \beta_2)}{2\ell_c^2}\overline{\zeta_1}.
\end{align*}

\paragraph{Terms of order $\cO(B^2\bB)$ and $\cO(B \bB^2)$.} Once again, by symmetry we have that $\Psi_{0,0,1,2,0}=-\mS_1\Psi_{0,0,2,1,0}$ where $\Psi_{0,0,2,1,0}$ solves
\begin{equation*}
0 = \T \Psi_{0,0,2,1,0} -\K * \left(\zeta_1^2\overline{\zeta_1}\right), \text{ with } \Psi_{0,0,2,1,0} \in \ker \Q.
\end{equation*}
We note that
\begin{equation*}
\K *  \left(\zeta_1^2\overline{\zeta_1}\right) (x) = \left[ \kappa_{0,1} x^3 -2\kappa_{2,1}x -\kappa_{3,1}\right]\rme^{\rmi \ell_c x}.
\end{equation*}
Here we used the fact that $\kappa_{1,1}=0$. As a consequence, we look for solutions of the form
\begin{equation*}
\Psi_{0,0,2,1,0} (x) = \left(\delta_0+\delta_1 x +\delta_2 x^2 + \delta_3 x^3\right) x^2 e^{\rmi \ell_c x} +\psi_{0,0,2,1,0} (x), \text{ with } \psi_{0,0,2,1,0} \in \E_0.
\end{equation*}
And we find the system satisfied by $(\delta_j)_{j=0,\cdots,3}$
\begin{align*}
10 \mu_c  \kappa_{2,1}\delta_3 -\kappa_{0,1}&=0,\\
 -10 \delta_3 \kappa_{3,1} +6 \delta_2 \kappa_{2,1} &=0,\\
 5 \mu_c \delta_3\kappa_{4,1} -4 \mu_c \delta_2 \kappa_{3,1}+(2-3\mu_c \delta_1) \kappa_{2,1} &=0,\\
 - \mu_c \delta_3\kappa_{5,1} + \mu_c \delta_2 \kappa_{4,1}+(1-\mu_c \delta_1) \kappa_{3,1}+\mu_c \delta_0 \kappa_{2,1} &=0,
\end{align*}
which can be solved recursively
\begin{equation*}
\delta_3=\frac{1}{10}\frac{\kappa_{0,1}^2}{\kappa_{2,1}}, \quad \delta_2 = \frac{1}{6} \frac{\kappa_{3,1}\kappa_{0,1}^2}{\kappa_{2,1}^2}, \quad \delta_1=\frac{1}{6}\frac{4\kappa_{0,1}\kappa_{2,1}^3-8\kappa_{3,1}^2\kappa_{0,1}^2+\kappa_{4,1} \kappa_{0,1}^2 \kappa_{2,1}}{\kappa_{2,1}^3},
\end{equation*}
and
\begin{equation*}
\delta_0 = \frac{1}{10}\frac{\kappa_{5,1}\kappa_{0,1}^2}{\kappa_{2,1}^2}-\frac{1}{6}\frac{\kappa_{3,1}\kappa_{0,1}\kappa_{2,1}^3+8\kappa_{3,1}^3\kappa_{2,1}^2}{\kappa_{2,1}^4}.
\end{equation*}
Once again, similar computations as above lead to
\begin{equation*}
\psi_{0,0,2,1,0}  = \frac{3(-\ell_c \delta_0+\rmi \delta_1)}{2\ell_c^3}\left( \zeta_0 - \overline{\zeta_0} \right) + \frac{4\rmi \ell_c \delta_0+3\delta_1}{2\ell_c^2}\zeta_1+ \frac{2\rmi \ell_c \delta_0+3\delta_1}{2\ell_c^2}\overline{\zeta_1}.
\end{equation*}

\paragraph{Terms of order $\cO(A^2B)$ and $\cO(\bA^2 \bB)$.} By symmetry we have that $\Psi_{0,2,0,1,0}=-\mS_1\Psi_{2,0,1,0,0}$ where $\Psi_{2,0,1,0,0}$ solves
\begin{equation*}
0 = \T \Psi_{2,0,1,0,0} -\K * \left(\zeta_0^2\zeta_1\right), \text{ with } \Psi_{2,0,1,0,0} \in \ker \Q,
\end{equation*}
and 
\begin{equation*}
\K * \left(\zeta_0^2\zeta_1\right)(x)= \left[ \kappa_{0,3} x - \kappa_{1,3}\right] \rme^{3\rmi\ell_cx}.
\end{equation*}
As a consequence, we look for solutions of the form
\begin{equation*}
\Psi_{2,0,1,0,0} (x) = \left(\gamma_0+\gamma_1 x \right) e^{3\rmi \ell_c x} +\psi_{2,0,1,0,0} (x), \text{ with } \psi_{2,0,1,0,0} \in \E_0.
\end{equation*}
Collecting terms of same order, we find a recursive system of equations
\begin{align*}
\gamma_1\left( -1+\mu_c\kappa_{0,3} \right) -\kappa_{0,3}&=0,\\
\gamma_0\left( -1+\mu_c \kappa_{0,3}\right) +\left(1-\mu_c \gamma_1\right) \kappa_{1,3} &=0.
\end{align*}
From which, we get
\begin{align*}
\gamma_1 &= \frac{\kappa_{0,3}}{ -1+\kappa_{0,3}/\kappa_{0,1}},\\
\gamma_0&=  \frac{ \kappa_{1,3}}{\left(-1+ \kappa_{0,3}/\kappa_{0,1} \right)^2}.
\end{align*}
One can also check that
\begin{equation*}
\psi_{2,0,1,0,0}=\left(\frac{4\ell_c \gamma_0-6\rmi \gamma_1}{\ell_c}\right)\zeta_0+\left(\frac{-5\ell_c \gamma_0+6\rmi \gamma_1}{\ell_c}\right)\overline{\zeta_0}-(8\gamma_1+8\rmi \ell_c\gamma_0)\zeta_1 -(5\gamma_1+4\rmi \ell_c\gamma_0)\overline{\zeta_1}.
\end{equation*}

\paragraph{Terms of order $\cO(A^2\bB)$ and $\cO(\bA^2 B)$.} By symmetry we have that $\Psi_{0,2,1,0,0}=-\mS_1\Psi_{2,0,0,1,0}$ where $\Psi_{2,0,0,1,0}\in \ker \Q$ solves
\begin{equation*}
0 = \T\Psi_{2,0,0,1,0} -\K * \left(\zeta_0^2\overline{\zeta_1}\right),
\end{equation*}
and $\K * \left(\zeta_0^2\overline{\zeta_1}\right) = \kappa_{0,1} \zeta_1$, so that 
\begin{equation*}
\Psi_{2,0,0,1,0}(x) = (\widetilde{\alpha_2} x^2 +\widetilde{\alpha_1} x) \zeta_1(x) + \psi_{2,0,0,1,0}(x),
\end{equation*}
for some $\psi_{2,0,0,1,0}\in \E_0$, and we get
\begin{align*}
-2\mu_c\widetilde{\alpha_2} \kappa_{2,1}-\kappa_{0,1}&=0,\\
\mu_c\widetilde{\alpha_1} \kappa_{2,1} - \mu_c \widetilde{\alpha_2} \kappa_{3,1}=0.
\end{align*}
From which, we deduce
\begin{equation*}
\widetilde{\alpha_2} =- \frac{\kappa_{0,1}^2}{2 \kappa_{2,1}},\quad \text{ and } \quad \widetilde{\alpha_1} = - \frac{ \kappa_{3,1} \kappa_{0,1}^2}{2\kappa_{2,1}^2}.
\end{equation*}
One also gets
\begin{equation*}
\psi_{2,0,0,1,0} = \frac{3\rmi\widetilde{\alpha_2}-3\ell_c \widetilde{\alpha_1} }{2\ell_c^3} \left(\zeta_0 -\overline{\zeta_0} \right)+\frac{4\rmi \ell_c\widetilde{\alpha_1}  +3 \widetilde{\alpha_2}}{2\ell_c^2}\zeta_1 
 + \frac{ 2\rmi \ell_c\widetilde{\alpha_1}  +3 \widetilde{\alpha_2}}{2\ell_c^2}\overline{\zeta_1}.
\end{equation*}

\paragraph{Terms of order $\cO(A\bA B)$ and $\cO(A\bA \bB)$.} By symmetry we have that $\Psi_{1,1,0,1,0}=-\mS_1\Psi_{1,1,1,0,0}$ where $\Psi_{1,1,1,0,0}\in \ker\Q$ solves
\begin{equation*}
0 = \T\Psi_{1,1,1,0,0} -2\K * \left(\zeta_1\right),
\end{equation*}
and $\K * \left(\zeta_1\right) = \kappa_{0,1} \zeta_1$. Using the previous computations, we find that 
\begin{equation*}
\Psi_{1,1,1,0,0}(x) = 2(\widetilde{\alpha_2} x^2 +\widetilde{\alpha_1} x) \zeta_1(x) + \psi_{1,1,1,0,0}(x), \text{ with } \psi_{1,1,1,0,0} \in \E_0
\end{equation*}
where
\begin{equation*}
\psi_{1,1,1,0,0} =2\psi_{2,0,0,1,0}.
\end{equation*}

\paragraph{Terms of order $\cO(A B \bB)$ and $\cO(\bA B \bB)$.} By symmetry we have that $\Psi_{0,1,1,1,0}=\mS_1\Psi_{1,0,1,1,0}$ where $\Psi_{1,0,1,1,0}$ solves
\begin{equation*}
0 = \T\Psi_{1,0,1,1,0} -2\K * \left(\zeta_0 \zeta_1 \overline{\zeta_1}\right), \text{ with } \Psi_{1,0,1,1,0} \in \ker \Q,
\end{equation*}
and 
\begin{equation*}
\K * \left(\zeta_0 \zeta_1 \overline{\zeta_1}\right)(x) = \left[ \kappa_{0,1} x^2 +\kappa_{2,1}\right] \rme^{\rmi \ell_c x},
\end{equation*}
so that 
\begin{equation*}
\Psi_{1,0,1,1,0}(x) = (\omega_2 x^2 +\omega_1 x +\omega_0) x^2\zeta_0(x) + \psi_{1,0,1,1,0}(x),
\end{equation*}
with $\psi_{1,0,1,1,0} \in \E_0$ and where $(\omega_j)_{j=0,1,2}$ solves
\begin{align*}
6\mu_c \omega_2 \kappa_{2,1} - 2\kappa_{0,1}&=0,\\
-4\mu_c \omega_2 \kappa_{3,1}-3\mu_c\kappa_{2,1}\omega_1&=0,\\
\mu_c\omega_2 \kappa_{4,1}-\mu_c \omega_1 \kappa_{3,1}+(\mu_c\omega_0-2) \kappa_{2,1}&=0.
\end{align*}
As a consequence, we have
\begin{equation*}
\omega_2 =\frac{1}{3} \frac{\kappa_{0,1}^2}{\kappa_{2,1}}, \quad \omega_1 =  -\frac{4}{9}\frac{\kappa_{0,1}^2\kappa_{3,1}}{\kappa_{2,1}^2}, \text{ and } \omega_0=2\kappa_{0,1}-\frac{4}{9}\frac{\kappa_{0,1}^2\kappa_{3,1}^2}{\kappa_{2,1}^3}-\frac{1}{3}\frac{\kappa_{4,1}\kappa_{0,1}^2}{\kappa_{2,1}^2}.
\end{equation*}
Finally, one can compute $\psi_{1,0,1,1,0}$ and we have
\begin{equation*}
\psi_{1,0,1,1,0}= \frac{3\rmi(\omega_1+\rmi \ell_c \omega_0)}{2\ell_c^3}(\zeta_0-\overline{\zeta_0})+\frac{3\omega_1+4\rmi \ell_c \omega_0}{2\ell_c^2}\zeta_1+\frac{3\omega_1+2\rmi \ell_c \omega_0}{2\ell_c^2}\overline{\zeta_1}.
\end{equation*}

\paragraph{Terms of order $\cO(\bA B^2)$ and $\cO(A \bB^2)$.} 
By symmetry we have that $\Psi_{1,0,0,2,0}=\mS_1\Psi_{0,1,2,0,0}$ where $\Psi_{0,1,2,0,0}\in \ker\Q$ solves
\begin{equation*}
0 = \T\Psi_{0,1,2,0,0} -\K * \left(\overline{\zeta_0}\zeta_1^2\right),
\end{equation*}
and 
\begin{equation*}
\K * \left(\overline{\zeta_0}\zeta_1^2\right)(x) = \left[ \kappa_{0,1} x^2 +\kappa_{2,1}\right] \rme^{\rmi \ell_c x},
\end{equation*}
so that $\Psi_{0,1,2,0,0}=\Psi_{1,0,1,1,0}$.

\paragraph{Terms of order $\cO(A B^2)$ and $\cO(\bA \bB^2)$.} 
By symmetry we have that $\Psi_{0,1,0,2,0}=\mS_1\Psi_{1,0,2,0,0}$ where $\Psi_{1,0,2,0,0}$ solves
\begin{equation*}
0 = \T\Psi_{1,0,2,0,0} -\K * \left(\zeta_0\zeta_1^2\right), \text{ with } \in \Psi_{1,0,2,0,0} \in \ker\Q,
\end{equation*}
and 
\begin{equation*}
\K * \left(\zeta_0\zeta_1^2\right)(x) = \left[ \kappa_{0,3} x^2-2\kappa_{1,3}x +\kappa_{2,3}\right] \rme^{3\rmi \ell_c x},
\end{equation*}
so that 
\begin{equation*}
\Psi_{1,0,2,0,0}(x) = (\rho_2 x^2 +\rho_1 x +\rho_0)\rme^{3\rmi \ell_c x}  +\psi_{1,0,2,0,0}(x) ,
\end{equation*}
where $\psi_{1,0,2,0,0}\in \E_0$ and 
\begin{align*}
\rho_2\left( -1+\mu_c\kappa_{0,3}\right)-\kappa_{0,3}&=0,\\
\rho_1\left( -1+\mu_c\kappa_{0,3}\right)-2\mu_c\rho_2\kappa_{1,3}+2\kappa_{1,3}&=0,\\
\rho_0\left( -1+\mu_c\kappa_{0,3}\right)+\mu_c\rho_2\kappa_{2,3}-\mu_c\rho_1\kappa_{1,3}-\kappa_{2,3}&=0.
\end{align*}
As a consequence, we get
\begin{equation*}
\rho_2=\frac{\kappa_{0,3}}{-1+\kappa_{0,3}/\kappa_{0,1}}, \quad \rho_1=\frac{2\kappa_{1,3}}{(-1+\kappa_{0,3}/\kappa_{0,1})^2}, \text{ and } \rho_0=-\frac{\kappa_{2,3}}{(-1+\kappa_{0,3}/\kappa_{0,1})^2}+\frac{2\kappa_{1,3}^2/\kappa_{0,1}}{(-1+\kappa_{0,3}/\kappa_{0,1})^3},
\end{equation*}
together with
\begin{align*}
\psi_{1,0,2,0,0}& = - \frac{-8\ell_c^2\rho_0+9\rho_2+12\rmi \ell_c \rho_1}{\ell_c^2}\zeta_0 + \frac{-10\ell_c^2\rho_0+9\rho_2+12\rmi \ell_c \rho_1}{\ell_c^2}\overline{\zeta_0} \\
&\qquad - \frac{8\rmi\ell_c^2 \rho_0-5\rmi \rho_2+8\ell_c\rho_1}{\ell_c}\zeta_1- \frac{4\rmi\ell_c^2 \rho_0-4\rmi \rho_2+5\ell_c\rho_1}{\ell_c}\overline{\zeta_1}.
\end{align*}

\subsection{Computations of reduced vector field at order 3}

We compute $\frac{\md}{\md x} \Q\left(\Psi(u_0(\cdot +x))\right)|_{x=0}$, the reduced vector field, induced by the reduction function $\Psi$.

We have that
\begin{align*}
\frac{\md}{\md x} \Q \left( (\cdot+x)^2\rme^{\rmi \ell_c (\cdot+x)}\right)|_{x=0}&=\left(\frac{3\rmi}{2\ell_c},-\frac{3\rmi}{2\ell_c},4,1  \right),\\
\frac{\md}{\md x} \Q \left( (\cdot+x)^3\rme^{\rmi \ell_c (\cdot+x)}\right)|_{x=0}&=\left(\frac{6}{\ell_c^2},-\frac{6}{\ell_c^2},-\frac{15\rmi}{2\ell_c},-\frac{9\rmi}{2\ell_c}  \right),\\
\frac{\md}{\md x} \Q \left( (\cdot+x)^4\rme^{\rmi \ell_c (\cdot+x)}\right)|_{x=0}&=\left(-\frac{6\rmi}{\ell_c^3},\frac{6\rmi}{\ell_c^3},-\frac{6}{\ell_c^2},-\frac{6}{\ell_c^2}  \right),\\
\frac{\md}{\md x} \Q \left( (\cdot+x)^5\rme^{\rmi \ell_c (\cdot+x)}\right)|_{x=0}&=(0,0,0,0),
\end{align*}
with
\begin{align*}
\frac{\md}{\md x} \Q \left( (\cdot+x)^2\rme^{-\rmi \ell_c (\cdot+x)}\right)|_{x=0}&=\left(\frac{3\rmi}{2\ell_c},-\frac{3\rmi}{2\ell_c},1,4  \right),\\
\frac{\md}{\md x} \Q \left( (\cdot+x)^3\rme^{-\rmi \ell_c (\cdot+x)}\right)|_{x=0}&=\left(-\frac{6}{\ell_c^2},\frac{6}{\ell_c^2},\frac{9\rmi}{2\ell_c},\frac{15\rmi}{2\ell_c}  \right),\\
\frac{\md}{\md x} \Q \left( (\cdot+x)^4\rme^{-\rmi \ell_c (\cdot+x)}\right)|_{x=0}&=\left(-\frac{6\rmi}{\ell_c^3},\frac{6\rmi}{\ell_c^3},-\frac{6}{\ell_c^2},-\frac{6}{\ell_c^2}  \right),\\
\frac{\md}{\md x} \Q \left( (\cdot+x)^5\rme^{-\rmi \ell_c (\cdot+x)}\right)|_{x=0}&=(0,0,0,0).
\end{align*}
Furthermore, similar computations lead to
\begin{align*}
\frac{\md}{\md x} \Q \left( \rme^{3\rmi \ell_c (\cdot+x)}\right)|_{x=0}&=\left(-12\rmi \ell_c, 15\rmi \ell_c, -24 \ell_c^2,-12 \ell_c^2  \right),\\
\frac{\md}{\md x} \Q \left( (\cdot+x)\rme^{3\rmi \ell_c (\cdot+x)}\right)|_{x=0}&=\left(-22,23,32\rmi \ell_c,19\rmi \ell_c \right),\\
\frac{\md}{\md x} \Q \left( (\cdot+x)^2\rme^{3\rmi \ell_c (\cdot+x)}\right)|_{x=0}&=\left(\frac{51\rmi}{2\ell_c},-\frac{51\rmi}{2\ell_c},31,22  \right),\\
\frac{\md}{\md x} \Q \left( (\cdot+x)^3\rme^{3\rmi \ell_c (\cdot+x)}\right)|_{x=0}&=\left(\frac{18}{\ell_c^2},-\frac{18}{\ell_c^2},-\frac{39\rmi}{2\ell_c},\frac{33\rmi}{2\ell_c}\right).
\end{align*}
with
\begin{align*}
\frac{\md}{\md x} \Q \left( \rme^{-3\rmi \ell_c (\cdot+x)}\right)|_{x=0}&=\left(-15\rmi \ell_c, 12\rmi \ell_c, -12 \ell_c^2,-24\ell_c^2  \right),\\
\frac{\md}{\md x} \Q \left( (\cdot+x)\rme^{-3\rmi \ell_c (\cdot+x)}\right)|_{x=0}&=\left(23,-22,-19\rmi \ell_c,-32\rmi \ell_c, \right),\\
\frac{\md}{\md x} \Q \left( (\cdot+x)^2\rme^{-3\rmi \ell_c (\cdot+x)}\right)|_{x=0}&=\left(\frac{51\rmi}{2\ell_c},-\frac{51\rmi}{2\ell_c},22,31  \right),\\
\frac{\md}{\md x} \Q \left( (\cdot+x)^3\rme^{-3\rmi \ell_c (\cdot+x)}\right)|_{x=0}&=\left((-\frac{18}{\ell_c^2},\frac{18}{\ell_c^2},\frac{33\rmi}{2\ell_c},\frac{39\rmi}{2\ell_c}\right).
\end{align*}

\bibliography{plain}

\end{document}